\documentclass{article}
\usepackage{verbatim}
\usepackage{amsmath}
\usepackage[mathscr]{eucal}
\usepackage{amsthm}
\usepackage{amssymb}
\usepackage{amsxtra}
\usepackage[english,german]{babel}
\theoremstyle{plain}
\newtheorem{Theo}{Theorem}[section]
\newtheorem{Prop}[Theo]{Proposition}
\newtheorem{Lem}[Theo]{Lemma}
\newtheorem{Cor}[Theo]{Corollary}
\newtheorem{Main}{Theorem}
\newtheorem{MainCor}{Corollary}
\theoremstyle{definition}
\newtheorem{Def}[Theo]{Definition}
\theoremstyle{remark}

\newcommand{\gl}{\mathfrak{g}\mathfrak{l}}
\newcommand{\g}{\mathfrak{g}}
\newcommand{\rank}{\operatorname{rank}}
\makeatletter
	
	\@addtoreset{equation}{section}
\makeatother
\begin{document}

\selectlanguage{english}
\title{Schur-Weyl reciprocity for the $q$-analogue of the alternating group}
\author{Hideo Mitsuhashi}

\date{}

\maketitle

\begin{center}
Department of Information Technology \\
Kanagawa Prefectural Junior College for Industrial Technology \\
2--4--1 Nakao, Asahi--ku, Yokohama--shi, Kanagawa--ken 241--0815, Japan
\end{center}

\begin{abstract}
In this paper, we establish Schur-Weyl reciprocity for the $q$-analogue of the alternating 
group. 
We analyze the sign $q$-permutation representation of the Hecke algebra 
$\mathcal{H}_{\mathbb{Q}(q),r}(q)$ on the $r$th tensor product of $\mathbb{Z}_2$-graded 
$\mathbb{Q}$-vector space $V=V_{\overline{0}}{\otimes}V_{\overline{1}}$ in detail, 
and examine its restriction to the $q$-analogue of the alternating group 
$\mathcal{H}_{\mathbb{Q}(q),r}^1(q)$. 
In consequence, we find out that if ${\dim}V_{\overline{0}}={\dim}V_{\overline{1}}$, 
then the centralizer of $\mathcal{H}_{\mathbb{Q}(q),r}^1(q)$ is a $\mathbb{Z}_2$-crossed product of 
the centralizer of $\mathcal{H}_{\mathbb{Q}(q),r}(q)$ and obtain Schur-Weyl reciprocity between 
$\mathcal{H}_{\mathbb{Q}(q),r}^1(q)$ and its centralizer. 
Though the structure of the centralizer is more complicated for the case 
${\dim}V_{\overline{0}}{\neq}{\dim}V_{\overline{1}}$, 
we obtain some results about the case. 
When $q=1$, Regev has proved Schur-Weyl reciprocity for alternating groups in \cite{Regev}. 
Therefore, our result can be regarded as an extension of Regev's work.
\end{abstract}

\renewcommand{\theenumi}{\arabic{enumi}}
\renewcommand{\labelenumi}{(\theenumi)}
\pagestyle{plain}
\markboth{}{}
\pagenumbering{arabic}

\section{Introduction}
The purpose of this study is to research Schur-Weyl reciprocity 
for the $q$-analogue of the alternating group. 
In our previous paper \cite{Mit2}, we established Schur-Weyl reciprocity between 
the Hecke algebra $\mathcal{H}_{\mathbb{Q}(q),r}(q)$ and the quantum super Lie algebra 
$U^\sigma_q\big{(}\gl(m,n)\big{)}$. 
In that paper, we defined the $q$-permutation representation of $\mathcal{H}_{\mathbb{Q}(q),r}(q)$, 
and showed that the image of the $q$-permutation representation is the centralizer of the image of 
the vector representation of the quantum super Lie algebra 
$U^\sigma_q\big{(}\gl(m,n)\big{)}$ on the $r$th tensor product of a $\mathbb{Z}_2$-graded 
$(m+n)$-dimensional $\mathbb{Q}(q)$-vector space $V=V_{\overline{0}}{\oplus}V_{\overline{1}}$. 
In this paper, we find out the centralizer of the $q$-analogue of the alternating group 
as the restriction of the $q$-permutation representation. 
When $q=1$, Regev has already shown Schur-Weyl reciprocity for the alternating group in \cite{Regev}. 
Hence our result is regarded as an extension of Regev's work. \par
Let $q$ be an indeterminate and $K=\mathbb{Q}(q)$. 
Let $(\pi_r,V^{{\otimes}r})$ be the $q$-permutation representation of $\mathcal{H}_{K,r}(q)$ 
(definition of the $q$-permutation representation is at (4.1)) and 
$(\rho_r,V^{{\otimes}r})$ the vector representation of $U^\sigma_q\big{(}\gl(m,n)\big{)}$ 
(definition of the vector representation is at (4.4)). 
We have proved in \cite{Mit2} that $\mathcal{A}_q=\pi_r\big{(}\mathcal{H}_{K,r}(q)\big{)}$ 
and $\mathcal{B}_q=\rho_r\big{(}U^\sigma_q\big{(}\gl(m,n)\big{)}\big{)}$ are full centralizers of each 
other, 
namely:
\begin{equation}
\mathcal{B}_q=\operatorname{End}_{\mathcal{A}_q}V^{{\otimes}r}{\quad}{\text{and}}{\quad}
\mathcal{A}_q=\operatorname{End}_{\mathcal{B}_q}V^{{\otimes}r}.
\end{equation}
Let $R_0$ be a commutative domain which includes an invertible element $q$. We further assume that 
$2$ and $q+q^{-1}$ are invertible elements of $R_0$. Then we can define the $q$-analogue of the 
alternating group $\mathcal{H}_{R_0,r}^1(q)$ (see Definition 3.1 and Proposition 3.6) in 
$\mathcal{H}_{R_0,r}(q)$. 
In \cite{Mit}, we defined the $q$-analogue of the alternating group as a subalgebra of 
Iwahori-Hecke algebra of type $A$ and obtained defining relations (see Proposition 3.7) 
for the first time. 
In this paper, we show that $\mathcal{H}_{R_0,r}(q)$ is isomorphic to the $\mathbb{Z}_2$-crossed 
product which is obtained from the crossed system 
$(\mathcal{H}_{R_0,r}^1(q),\mathbb{Z}_2,\psi_0,\alpha_0)$ 
(definition of $\psi_0$ and $\alpha_0$ are at (3.2) and (3.3) respectively). 
\renewcommand{\theMain}{3.9}
\begin{Main}
$\mathcal{H}_{R_0,r}(q)$ is isomorphic to the $\mathbb{Z}_2$-crossed product 
$\mathcal{H}_{R_0,r}^1(q)_{\alpha_0}^{\psi_0}[\mathbb{Z}_2]$ as $R_0$-algebras. 
\end{Main}
Let $\mathcal{C}_q=\pi_r\big{(}\mathcal{H}_{K,r}^1(q)\big{)}$ and 
$\mathcal{D}_q=\operatorname{End}_{\mathcal{C}_q}V^{{\otimes}r}$. 
Our main subject is to solve the relation between $\mathcal{D}_q$ and $\mathcal{B}_q$. 
From (1.1), one can immediately see that  $\mathcal{B}_q{\subseteq}\mathcal{D}_q$. 
But the structure of $\mathcal{D}_q$ is not trivial. 
Indeed, the structure of $\mathcal{D}_q$ depends on the dimensions $m=\dim_{K}V_{\overline{0}}$ 
and $n=\dim_{K}V_{\overline{1}}$. 
In this paper, we show that if $m=n$, then 
$\mathcal{D}_q$ is isomorphic to 
$\mathbb{Z}_2$-crossed product which is obtained from the crossed system 
$(\mathcal{B}_q,\mathbb{Z}_2,\psi_1,\alpha_1)$ 
(definition of $\psi_0$ and $\alpha_0$ are at (5.4) and (5.5) respectively). 
\renewcommand{\theMain}{5.6}
\begin{Main}
If $m=n$, then $\mathcal{D}_q$ is isomorphic to the $\mathbb{Z}_2$-crossed product 
${\mathcal{B}_q}_{\alpha_1}^{\psi_1}[\mathbb{Z}_2]$ as $K$-algebras. 
\end{Main}
From this theorem, we immediately obtain that $\dim_{K}\mathcal{D}_q=2\dim_{K}\mathcal{B}_q$. 
Moreover, we show Schur-Weyl reciprocity for $\mathcal{H}_{K,r}^1(q)$. \par
\renewcommand{\theMain}{5.8}
\begin{Main}
$\operatorname{End}_{\mathcal{C}_q}V^{{\otimes}r}
=\mathcal{D}_q$ and $\operatorname{End}_{\mathcal{D}_q}V^{{\otimes}r}
=\mathcal{C}_q$ hold. 
\end{Main}
In the general case, the matter is more complicated, but we can find out to some extent if we 
exchange the base field from $K$ to its algebraic closure $\bar{K}$. 
Let $\bar{U}^\sigma_q\big{(}\gl(m,n)\big{)}=U^\sigma_q\big{(}\gl(m,n)\big{)}
{\otimes_{K}}\bar{K}$, 
$\bar{\mathcal{A}}_q=\mathcal{A}_q{\otimes_{K}}\bar{K}$, 
$\bar{\mathcal{B}}_q=\mathcal{B}_q{\otimes_{K}}\bar{K}$ and 
$\bar{\mathcal{C}}_q=\mathcal{C}_q{\otimes_{K}}\bar{K}$. 
Then we have the following theorem by the theory of semisimple algebras. 
\renewcommand{\theMain}{6.1}
\begin{Main}
$\bar{\mathcal{A}}_q$ and $\bar{\mathcal{C}}_q$ have direct sum decompositions 
$\bar{\mathcal{A}}_q=\bar{\mathcal{A}}_q^0{\oplus}\bar{\mathcal{A}}_q^1$ and 
$\bar{\mathcal{C}}_q=\bar{\mathcal{C}}_q^0{\oplus}\bar{\mathcal{C}}_q^1$ respectively, 
which are satisfy the following relations. 
\begin{enumerate}
\item
$\bar{\mathcal{C}}_q^0{\subseteq}\bar{\mathcal{A}}_q^0$ and 
$\dim_{\bar{K}}\bar{\mathcal{A}}_q^0
=2\dim_{\bar{K}}\bar{\mathcal{C}}_q^0$
\item
$\bar{\mathcal{C}}_q^1=\bar{\mathcal{A}}_q^1$. Especially 
$\dim_{\bar{K}}\bar{\mathcal{A}}_q^1
=\dim_{\bar{K}}\bar{\mathcal{C}}_q^1$. 
\end{enumerate}
\end{Main}
As a corollary, we obtain an anomalous phenomenon for non-super case as follows. 
\renewcommand{\theMainCor}{6.2}
\begin{MainCor}
Let $n=0$. If $m^2<r$, then $\bar{\mathcal{A}}_q=\bar{\mathcal{C}}_q$ and 
$\operatorname{End}_{\bar{\mathcal{C}}_q}\bar{V}^{{\otimes}r}=
\operatorname{End}_{\bar{\mathcal{A}}_q}\bar{V}^{{\otimes}r}$.
\end{MainCor}
The details of the decompositions of $\bar{\mathcal{A}}_q$ and $\bar{\mathcal{C}}_q$ are 
described in section 6. 
We also obtain the similar result to Theorem 6.1 about the endomorphism algebras 
$\operatorname{End}_{\bar{\mathcal{A}}_q}\bar{V}^{{\otimes}r}$ and 
$\operatorname{End}_{\bar{\mathcal{C}}_q}\bar{V}^{{\otimes}r}$; 
there exist two $\mathcal{H}_{\bar{K},r}(q)$-submodules $W_0$ and $W_1$ which 
satisfy the following properties. 
\renewcommand{\theMainCor}{6.3}
\begin{MainCor}
$\operatorname{End}_{\bar{\mathcal{C}}_q}W_0{\supseteq}\operatorname{End}_{\bar{\mathcal{A}}_q}W_0$ and 
$\operatorname{End}_{\bar{\mathcal{C}}_q}W_1=\operatorname{End}_{\bar{\mathcal{A}}_q}W_1$.
\end{MainCor}
Although the relation between $\bar{\mathcal{A}}_q$ and $\bar{\mathcal{C}}_q$ is made clear by 
(6.4)-(6.7), that between $\operatorname{End}_{\bar{\mathcal{A}}_q}W$ and 
$\operatorname{End}_{\bar{\mathcal{C}}_q}W$ is not clear except for Corollary 6.2 at this point. 

\section{Preliminaries}
Let $R$ be a commutative ring with $1$ and $G$ a group. 
In this section, we shall review the definition and some properties about $G$-crossed products. 
A full account about $G$-graded algebras and $G$-crossed products is given in \cite{N-O}.
\begin{Def}[$G$-graded algebra]
An $R$-algebra $A$ is said to be \it $G$-graded \rm if there exist a family of $R$-submodules 
$\{A_{\sigma}|\sigma{\in}G\}$ of $A$ indexed by elements of $G$ which satisfies 
the following two conditions: 
\renewcommand{\theenumi}{\arabic{enumi}}
\renewcommand{\labelenumi}{(G\theenumi)}
\begin{enumerate}
\item
$A=\bigoplus_{\sigma{\in}G}A_{\sigma}$,
\item
$A_{\sigma}A_{\tau}{\subseteq}A_{\sigma\tau}$ for $\sigma,\tau{\in}G$.
\end{enumerate}
\renewcommand{\theenumi}{\arabic{enumi}}
\renewcommand{\labelenumi}{(\theenumi)}
Moreover, $A$ is said to be \it strongly $G$-graded \rm when (G2) is replaced by 
the following condition: 
\renewcommand{\theenumi}{\arabic{enumi}}
\renewcommand{\labelenumi}{(G\theenumi)}
\begin{enumerate}
\item[(G'2)]
$A_{\sigma}A_{\tau}=A_{\sigma\tau}$ for $\sigma,\tau{\in}G$.
\end{enumerate}
\renewcommand{\theenumi}{\arabic{enumi}}
\renewcommand{\labelenumi}{(\theenumi)}
\end{Def}
We notice that if $A=\bigoplus_{\sigma{\in}G}A_{\sigma}$ is a $G$-graded algebra, then $A_{1_G}$ 
($1_G$ means the identity element of $G$) is a subalgebra of $A$ and $1_A{\in}A_{1_G}$. 
\begin{Def}[$G$-crossed product]
A $G$-graded $R$-algebra $A=\oplus_{\sigma{\in}G}A_{\sigma}$ is said to be a 
$G$-crossed product if each $A_{\sigma}$ has an invertible element. 
\end{Def}
We notice that a $G$-crossed product is a strongly $G$-graded algebra. 
Indeed, if $A$ is a $G$-crossed product, then for an invertible element 
$u_{\sigma}{\in}A_{\sigma}$, $u_{\sigma}^{-1}{\in}A_{\sigma^{-1}}$ and 
$1_A=u_{\sigma}u_{\sigma}^{-1}{\in}A_{\sigma}A_{\sigma^{-1}}$. So, we have 
$A_{\sigma\tau}=1_AA_{\sigma\tau}{\subseteq}(A_{\sigma}A_{\sigma^{-1}})A_{\sigma\tau}
=A_{\sigma}(A_{\sigma^{-1}}A_{\sigma\tau}){\subseteq}A_{\sigma}A_{\tau}$. 
\begin{Def}[crossed system]
Let $A$ be an algebra and $G$ a group. 
Suppose that there exist two maps 
\begin{equation*}
\psi:G{\longrightarrow}\operatorname{Aut}(A),
\end{equation*}
(we denote $\psi(\sigma)(a)$ by ${^{\sigma}a}$ for brevity), and 
\begin{equation*}
\alpha:G{\times}G{\longrightarrow}A^{\times},
\end{equation*}
where $A^{\times}$ is the multiplicative group of units of $A_1$, 
which satisfy the relations: 
\begin{eqnarray}
{^\sigma}({^\tau}a)
&=&\alpha(\sigma,\tau){^{({\sigma}{\tau})}}a\alpha(\sigma,\tau)^{-1} \\
{^{\sigma_1}}\alpha(\sigma_2,\sigma_3)\alpha(\sigma_1,\sigma_2\sigma_3)
&=&\alpha(\sigma_1,\sigma_2)\alpha(\sigma_1\sigma_2,\sigma_3) \\
\alpha(\sigma,1)
&=&\alpha(1,\sigma)=1, 
\end{eqnarray}
for $\sigma,\tau,\sigma_1,\sigma_2,\sigma_3{\in}G,a{\in}A$. 
$(A,G,\psi,\alpha)$ is said to be a \it crossed system\rm. 
$\psi$ is called a \it weak action \rm of $G$ on $A$, and $\alpha$ is called a 
$\psi$-\it cocycle\rm. 
\end{Def}
We denote by $A_{\alpha}^{\psi}[G]$ the free left $A$-module with the basis 
$\{u_{\sigma}|\sigma{\in}G\}$ and the following multiplication: 
\begin{equation}
(a_1u_{\sigma})(a_2u_{\tau})=a_1{^\sigma}a_2\alpha(\sigma,\tau)u_{\sigma\tau},
\end{equation}
for $a_1,a_2{\in}A,\sigma,\tau{\in}G$. 
\begin{Prop}[\cite{N-O}Proposition 1.4.1]
$A_{\alpha}^{\psi}[G]$ is a $G$-crossed product. 
\end{Prop}
\begin{Prop}[\cite{N-O}Proposition 1.4.2]
Every $G$-crossed product is of the form $A_{\alpha}^{\psi}[G]$ for some algebra $A$, 
some weak action $\psi$ and some $\psi$-cocycle $\alpha$. 
\end{Prop}
When $G$ is finite and a strongly $G$-graded algebra $A=\bigoplus_{\sigma{\in}G}A_{\sigma}$ 
is finitely generated over $R$ as modules, 
$A$ is said to have a \it $G$-graded Clifford system \rm $\{A_{\sigma}|\sigma{\in}G\}$ if 
$A$ satisfies (C1). 
\renewcommand{\theenumi}{\arabic{enumi}}
\renewcommand{\labelenumi}{(C\theenumi)}
\begin{enumerate}
\item
For each $\sigma{\in}G$, there exists an invertible element $a_{\sigma}{\in}A$ such that 
$A_{\sigma}=a_{\sigma}A_{1_G}=A_{1_G}a_{\sigma}$. 
\end{enumerate}
\renewcommand{\theenumi}{\arabic{enumi}}
\renewcommand{\labelenumi}{(\theenumi)}
It is clear that such $a_{\sigma}$ is in $A_{\sigma}$. 
An exposition about group graded Clifford systems can be found in \cite{C-R1}, section 11C. 

\section{The $q$-analogue of the alternating group and its representation}
Let $(W,S=\{s_1,\hdots,s_r\})$ be a Coxeter system of rank $r$. 
Let $R_0$ be a commutative domain with $1$, and let $q_i(i=1,\hdots,r)$ be any invertible 
elements of $R_0$ such that $q_i=q_j$ if $s_i$ is conjugate to $s_j$ in $W$. 
Further we assume that $2$ and $q_i+q_i^{-1}$($i=1,2,\hdots,r$) are invertible elements of $R_0$. 
The Iwahori-Hecke algebra 
$\mathcal{H}_{R_0}(W,S)$ is an $R_0$-algebra generated by $\{T_{s_i}|s_i{\in}S\}$ with the 
defining relations: 
\renewcommand{\theenumi}{\arabic{enumi}}
\renewcommand{\labelenumi}{(H\theenumi)}
\begin{enumerate}
\item
$T_{s_i}^2 = (q_i-q_i^{-1})T_{s_i}+1$ \qquad if $i=1,2,\hdots,r$,
\item
$(T_{s_i}T_{s_j})^{k_{ij}}=(T_{s_j}T_{s_i})^{k_{ij}}$ \qquad if $m_{ij}=2k_{ij}$,
\item
$(T_{s_i}T_{s_j})^{k_{ij}}T_{s_i}=(T_{s_j}T_{s_i})^{k_{ij}}T_{s_j}$ \qquad if $m_{ij}=2k_{ij}+1$,
\end{enumerate}
\renewcommand{\theenumi}{\arabic{enumi}}
\renewcommand{\labelenumi}{(\theenumi)}
where $m_{ij}$ is the order of $s_is_j$ in $W$. 
We write $T_i=T_{s_i}$ for brevity. \par
If $(W,S)$ is of type $A$ and of rank $r-1$, then $W$ is isomorphic to the symmetric group 
$\mathfrak{S}_r$. 
Furthermore, all the elements of $S$ are conjugate to each other, hence we may assume 
$q_1=\cdots=q_{r-1}=q$. 
The Iwahori-Hecke algebra $\mathcal{H}_{R_0,r}(q)=\mathcal{H}_{R_0}(W,S)$ of type $A$ has 
the defining relations: 
\renewcommand{\theenumi}{\arabic{enumi}}
\renewcommand{\labelenumi}{(A\theenumi)}
\begin{enumerate}
\item
$T_i^2 = (q-q^{-1})T_i+1$ \qquad if $i=1,2,\hdots,r-1$,
\item
$T_iT_{i+1}T_i = T_{i+1}T_iT_{i+1}$ \qquad if $i=1,2,\hdots,r-2$,
\item
$T_iT_j = T_jT_i$ \qquad if $|i-j|>1$.
\end{enumerate}
\renewcommand{\theenumi}{\arabic{enumi}}
\renewcommand{\labelenumi}{(\theenumi)}\par
Let \, $\hat{}$ \, be the Goldman involution. This is an involution on $\mathcal{H}_{R_0,r}(q)$ 
defined by 
\[
\hat{T_i}=(q-q^{-1})-T_i.
\]
\begin{Def}
We define $\mathcal{H}_{R_0,r}^{\pm1}(q)$ to be the eigenspaces of $\mathcal{H}_{R_0,r}(q)$ 
corresponding to the eigenvalues $\pm1$ of \, $\hat{}$ \, respectively. 
\end{Def}
We notice that $\mathcal{H}_{R_0,r}^1(q)$ is a subalgebra of $\mathcal{H}_{R_0,r}(q)$. 
Let $T'_i$ ($i=1,2,\hdots,r-1$) be the elements of $\mathcal{H}_{R_0,r}(q)$ defined by 
\[
T'_i=\dfrac{T_i-\hat{T}_i}{q+q^{-1}}= \dfrac{2T_i-(q-q^{-1})}{q+q^{-1}} 
\quad \text{for $i=1,2,\hdots,r-1$}.
\]
Then one can immediately check $\hat{T}'_i=-T'_i$. 
\begin{Prop}
$T'_i$($i=1,2,\hdots,r$) generate $\mathcal{H}_{R_0,r}(q)$ and satisfy the following 
defining relations: 
\renewcommand{\theenumi}{\arabic{enumi}}
\renewcommand{\labelenumi}{(A'\theenumi)}
\begin{enumerate}
\item
${T'_i}^2 = 1$ \qquad if $i=1,2,\hdots,r-1$,
\item
$T'_iT'_{i+1}T'_i = T'_{i+1}T'_iT'_{i+1}-\Big{(}\dfrac{q-q^{-1}}{q+q^{-1}}\Big{)}^2
(T'_i-T'_{i+1})$ \qquad if $i=1,2,\hdots,r-2$,
\item
$T'_iT'_j = T'_jT'_i$ \qquad if $|i-j|>1$.
\end{enumerate}
\renewcommand{\theenumi}{\arabic{enumi}}
\renewcommand{\labelenumi}{(\theenumi)}
\end{Prop}
\begin{proof}
From the equations 
\[
T_i=\dfrac{1}{2}\{(q+q^{-1})T'_i+(q-q^{-1})\} \quad \text{for $i=1,2,\hdots,r-1$},
\]
which are obtained from the definition of $T'_i$, 
we see that $T'_i$($i=1,2,\hdots,r$) generate $\mathcal{H}_{R_0,r}(q)$. 
The defining relations are obtained from a direct computation. 
\end{proof}
Consider the following sets of monomials: 
\begin{equation*}
\begin{split}
\mathcal{C}_1&=\{1,T_1\} \\
\mathcal{C}_2&=\{1,T_2,T_2T_1\} \\
\mathcal{C}_3&=\{1,T_3,T_3T_2,T_3T_2T_1\} \\
&: \\
\mathcal{C}_{r-1}&=\{1,T_{r-1},T_{r-1}T_{r-2},\hdots,T_{r-1}T_{r-2}{\cdots}T_1\}
\end{split}
\end{equation*}
We shall say that $M_1M_2{\cdots}M_{r-1}$ is a \it monomial in $T_i$-normal form \rm 
in $\mathcal{H}_{R_0,r}(q)$ if $M_i{\in}\mathcal{C}_i$ for $i=1,2,\hdots,r-1$. 
The following fact is well-known in the theory of the Iwahori-Hecke algebra. 
\begin{Prop}
$\operatorname{rank}_{R_0}\mathcal{H}_{R_0,r}(q)=r!$ and 
\[
\mathcal{H}_{R_0,r}(q)=\bigoplus_{M_i{\in}\mathcal{C}_i}R_0M_1M_2{\cdots}M_{r-1}
\]
\end{Prop}
We derive from this fact that all monomials in $T'_i$-normal form also constitute 
a basis of $\mathcal{H}_{R_0,r}(q)$. 
\begin{Prop}
Let $\mathcal{C}'_i=\{1,T'_i,T'_iT'_{i-1},\hdots,T'_iT'_{i-1}{\cdots}T'_1\}$ for $i=1,2,\hdots,r-1$. 
Then we have 
\[
\mathcal{H}_{R_0,r}(q)=\bigoplus_{M'_i{\in}\mathcal{C}'_i}R_0M'_1M'_2{\cdots}M'_{r-1}.
\]
\end{Prop}
\begin{proof}
Consider the map $f$ from $\{T_i\}_{i=1,2,\hdots,r-1}$ to $\{T'_i\}_{i=1,2,\hdots,r-1}$ which is 
defined by $f(T_i)=T'_i$. This map induces the $R_0$-endomorphism $\bar{f}$ of 
$\mathcal{H}_{R_0,r}(q)$. $\bar{f}$ is an $R_0$-isomorphism because the inverse $\bar{g}$ 
which is induced from the map $g$ from $\{T'_i\}_{i=1,2,\hdots,r-1}$ to 
$\{T_i\}_{i=1,2,\hdots,r-1}$ defined by 
\[
g(T'_i)=\dfrac{1}{2}\{(q+q^{-1})T_i+(q-q^{-1})\}
\]
exists. Hence we conclude that $M'_1M'_2{\cdots}M'_{r-1}$($M'_i{\in}\mathcal{C}'_i$), 
which are images of $M_1M_2{\cdots}M_{r-1}$($M_i{\in}\mathcal{C}_i$), are linearly independent and 
constitute a basis of $\mathcal{H}_{R_0,r}(q)$. 
\end{proof}
Let $\mathcal{E}_r$ (respectively $\mathcal{O}_r$) be the set of all monomials in 
$T'_i$-normal form in $\mathcal{H}_{R_0,r}(q)$ which are products of even (respectively odd) 
numbers of $T'_i$'s. 
Then the following holds. 
\begin{Lem}
$|\mathcal{E}_r|=|\mathcal{O}_r|=2^{-1}r!$ for $r>1$. 
\end{Lem}
\begin{proof}
The proof is done by induction on $r$. 
It is trivial for $r=2$. 
Let $M'_1M'_2{\cdots}M'_{r-1}{\in}\mathcal{E}_r$ with $r>2$ and $M'_i{\in}\mathcal{C}'_i$. 
Then $M'_1M'_2{\cdots}M'_{r-2}$ is considered as a monomial in $T'_i$-normal form in 
$\mathcal{H}_{R_0,r-1}(q)$. 
Let $(\mathcal{C}'_{r-1})^e$ (respectively $(\mathcal{C}'_{r-1})^o$) be the subset of 
$\mathcal{C}'_{r-1}$ which consists of products of even (respectively odd) numbers of $T'_i$'s. 
We can readily see that $M'_1M'_2{\cdots}M'_{r-2}{\in}\mathcal{E}_{r-1}$ if and only if 
$M'_{r-1}{\in}(\mathcal{C}'_{r-1})^e$. 
By induction, $|\mathcal{E}_{r-1}|=|\mathcal{O}_{r-1}|=2^{-1}(r-1)!$ and hence 
we obtain the following
\[
|\mathcal{E}_r|=|(\mathcal{C}'_{r-1})^e||\mathcal{E}_{r-1}|
+|(\mathcal{C}'_{r-1})^o||\mathcal{O}_{r-1}| 
=r{\times}2^{-1}(r-1)!=2^{-1}r!
\]
as desired. 
\end{proof}
Now we characterize $\mathcal{H}_{R_0,r}^1(q)$ as a $q$-analogue of the alternating group. 
\begin{Prop}
$\operatorname{rank}_{R_0}\mathcal{H}_{R_0,r}(q)=2\operatorname{rank}_{R_0}\mathcal{H}_{R_0,r}^1(q)$. 
Moreover $\mathcal{H}_{R_0,r}^1(q)$ is the subalgebra which consists of all the products of even 
numbers of $T'_i$'s. 
\end{Prop}
\begin{proof}
Let $\mathcal{H}_{R_0,r}^e(q)=\oplus_{M{\in}\mathcal{E}_r}R_0M$ and
$\mathcal{H}_{R_0,r}^o(q)=\oplus_{M{\in}\mathcal{O}_r}R_0M$. 
Then we can see immediately that 
$\mathcal{H}_{R_0,r}(q)=\mathcal{H}_{R_0,r}^e(q){\oplus}\mathcal{H}_{R_0,r}^o(q)$. 
Furthermore we obtain 
$\mathcal{H}_{R_0,r}^e(q)=\mathcal{H}_{R_0,r}^1(q)$ and 
$\mathcal{H}_{R_0,r}^o(q)=\mathcal{H}_{R_0,r}^{-1}(q)$ from the property $\hat{T}'_i=-T'_i$ . 
Combined with Lemma 3.5 we have 
$\operatorname{rank}_{R_0}\mathcal{H}_{R_0,r}(q)=2\operatorname{rank}_{R_0}\mathcal{H}_{R_0,r}^1(q)$. 
Let $\bar{\mathcal{H}}_{R_0,r}^1(q)$ be the set of all the linear combinations of products of 
even numbers of $T'_i$'s. 
Obviously $\mathcal{H}_{R_0,r}^1(q){\subseteq}\bar{\mathcal{H}}_{R_0,r}^1(q)$. 
From the defining relations (A'1)-(A'3), one can see that if a monomial in $T'_i$-normal form 
which consists of even (respectively odd) number of $T'_i$'s 
is expressed in a linear combination of other expressions, 
then each term consists of even (respectively odd) number of $T'_i$'s. 
Hence if we express an element of $\bar{\mathcal{H}}_{R_0,r}^1(q)$ by a linear combination 
of monomials in $T'_i$-normal form, each term is in $\mathcal{H}_{R_0,r}^1(q)$. 
Consequently $\bar{\mathcal{H}}_{R_0,r}^1(q)=\mathcal{H}_{R_0,r}^1(q)$. 
\end{proof}
When we suppose that $R_0=\mathbb{C}$ and take a limit $q{\rightarrow}1$, 
$\mathcal{H}_{\mathbb{C},r}^1(1)$ is isomorphic to the group algebra $\mathbb{C}[\mathfrak{A}_n]$ of 
the alternating group $\mathfrak{A}_n$. 
\begin{Theo}[\cite{Mit}]
$\mathcal{H}_{R_0,r}^1(q)$ is isomorphic to the $R_0$-algebra which is generated by $r-2$ elements 
$X_1,X_2,\hdots,X_{r-2}$ with the defining relations: 
\renewcommand{\theenumi}{\arabic{enumi}}
\renewcommand{\labelenumi}{(B\theenumi)}
\begin{enumerate}
\item
$X_1^3=-\Big{(}\dfrac{q-q^{-1}}{q+q^{-1}}\Big{)}^2(X_1^2-X_1)+1$,
\item
$X_i^2=1 \qquad \text{for $i>1$}$,
\item
$(X_{i-1}X_i)^3=-\Big{(}\dfrac{q-q^{-1}}{q+q^{-1}}\Big{)}^2\Big{\{}(X_{i-1}X_i)^2-X_{i-1}X_i\Big{\}}+1 
\qquad \text{for $i=2,3,\hdots ,r-2$}$,
\item
$(X_iX_j)^2=1 \qquad \text{whenever $|i-j|>1$}$.
\end{enumerate}
\renewcommand{\theenumi}{\arabic{enumi}}
\renewcommand{\labelenumi}{(\theenumi)}
An isomorphism is given by $X_i{\longrightarrow}T'_1T'_{i+1}$. 
\end{Theo}
Next, we shall show that $\mathcal{H}_{R_0,r}(q)$ is a $\mathbb{Z}_2$-crossed product. 
Since $T'_i$ has an inverse as itself, we can readily see from Proposition 3.6 that 
$\mathcal{H}_{R_0,r}^1(q)T'_1$ is an $R_0$-submodule of $\mathcal{H}_{R_0,r}(q)$ which 
consists of all linear combinations of products of odd 
numbers of $T'_i$'s. Therefore, we obtain a direct sum decomposition 
of left $\mathcal{H}_{R_0,r}^1(q)$-modules: 
\begin{equation}
\mathcal{H}_{R_0,r}(q)=\mathcal{H}_{R_0,r}^1(q){\oplus}\mathcal{H}_{R_0,r}^1(q)T'_1. 
\end{equation}
Let $\mathbb{Z}_2={\langle}1,-1{\rangle}$ be a multiplicative group. 
We define two maps $\psi_0$ and $\alpha_0$ to be 
\begin{equation}
\psi_0:\mathbb{Z}_2{\longrightarrow}\operatorname{Aut}(\mathcal{H}_{R_0,r}^1(q)),\quad 
\psi_0(1)(T)=T,\psi_0(-1)(T)=T'_1TT'_1,\quad T{\in}\mathcal{H}_{R_0,r}^1(q). 
\end{equation}
and
\begin{equation}
\alpha_0:\mathbb{Z}_2{\times}\mathbb{Z}_2{\longrightarrow}
(\mathcal{H}_{R_0,r}^1(q))^{\times},\quad 
\alpha_0(\sigma,\tau)=1 \quad\text{for all $\sigma,\tau{\in}\mathbb{Z}_2$}.
\end{equation}
Then we have the following immediately. 
\begin{Lem}
$\psi_0$ and $\alpha_0$ satisfy (2.1)--(2.3).
\end{Lem}
Thus, we obtain a $\mathbb{Z}_2$-crossed product 
$\mathcal{H}_{R_0,r}^1(q)_{\alpha_0}^{\psi_0}[\mathbb{Z}_2]$ 
from the crossed system $(\mathcal{H}_{R_0,r}^1(q),\mathbb{Z}_2,\psi_0,\alpha_0)$, 
where $\psi_0$ and $\alpha_0$ are given by (3.2) and (3.3) respectively. 
\begin{Theo}
$\mathcal{H}_{R_0,r}(q)$ is isomorphic to 
$\mathcal{H}_{R_0,r}^1(q)_{\alpha_0}^{\psi_0}[\mathbb{Z}_2]$ as $R_0$-algebras. 
\end{Theo}
\begin{proof}
Since both $\mathcal{H}_{R_0,r}(q)$ and 
$\mathcal{H}_{R_0,r}^1(q)_{\alpha_0}^{\psi_0}[\mathbb{Z}_2]$ are free left 
$\mathcal{H}_{R_0,r}^1(q)$-modules, we may define an isomorphism of 
$\mathcal{H}_{R_0,r}^1(q)$-modules 
\begin{equation*}
{\iota_0}:\mathcal{H}_{R_0,r}^1(q)_{\alpha_0}^{\psi_0}[\mathbb{Z}_2]{\longrightarrow}
\mathcal{H}_{R_0,r}(q),\quad 
{\iota_0}(u_1)=1,{\iota_0}(u_{-1})=T'_1
\end{equation*}
From (2.4) and (3.2) and (3.3), we can determine the multiplication law as follows. 
\begin{equation*}
(a_1u_{\sigma})(a_2u_{\tau})=
\begin{cases}
a_1a_2u_{\sigma\tau}& \text{if $\sigma=1$,}\\
a_1T'_1a_2T'_1u_{\sigma\tau}& \text{if $\sigma=-1$,}
\end{cases}
\end{equation*}
where $a_1,a_2{\in}\mathcal{H}_{R_0,r}^1(q)$ and $\sigma,\tau{\in}\mathbb{Z}_2$. 
Therefore, we get four formulas in $\mathcal{H}_{R_0,r}(q)$ 
\begin{equation*}
(a_11)(a_21)=a_1a_21,{\quad}(a_11)(a_2T'_1)=a_1a_2T'_1,{\quad}
(a_1T'_1)(a_21)=a_1T'_1a_2T'_1T'_1,{\quad}(a_1T'_1)(a_2T'_1)=a_1T'_1a_2T'_11,
\end{equation*}
which derive the conclusion that ${\iota_0}$ is an isomorphism of $R_0$-algebras. 
\end{proof}
We denote by $\bar{K}$ an algebraic closure of a field $K$. 
Let $q$ be an indeterminate and $K=\mathbb{Q}(q)$. 
We shall show (split) semisimplicity of $\mathcal{H}_{\bar{K},r}^1(q)$ and the branching rule 
from $\mathcal{H}_{\bar{K},r}(q)$ to $\mathcal{H}_{\bar{K},r}^1(q)$. 
The manner of proof given here is credited to K.Uno, who sent me a letter enclosing the outline of 
this proof. 
It is well-known that $\mathcal{H}_{\bar{K},r}(q)$ is split semisimple 
and that isomorphism classes of simple left $\mathcal{H}_{\bar{K},r}(q)$-modules are parametrized by 
Young diagrams of total size $n$. 
Let $\Lambda_r$ be the set of all Young diagrams of total size $r$. 
Then, provided that  $\{M_{q,\lambda}|\lambda\in\Lambda_r\}$ is a set of all isomorphism classes of 
simple left $\mathcal{H}_{\bar{K},r}(q)$-modules and that $d_{\lambda}={\deg}M_{q,\lambda}$, we may write 
\begin{equation}
\mathcal{H}_{\bar{K},r}(q)=\bigoplus_{\lambda{\in}\Lambda_r}I_{q,\lambda}{\quad}
\big{(}I_{q,\lambda}{\cong}\operatorname{Mat}_{d_{\lambda}}(\bar{K})\big{)},
\end{equation}
where each $I_{q,\lambda}$ is the homogeneous component corresponding to $\lambda$. 
$I_{q,\lambda}$ is isomorphic to a $d_{\lambda}{\times}d_{\lambda}$ matrix algebra 
$\operatorname{Mat}_{d_{\lambda}}(\bar{K})$ whose entries lie in $\bar{K}$. 
Since \, $\hat{}$ \, is an involution of $\mathcal{H}_{\bar{K},r}(q)$, for each $\lambda{\in}\Lambda_r$ 
there exists $\mu{\in}\Lambda_n$ such that 
$\hat{I}_{q,\lambda}=I_{q,\mu}$. 
Especially, $d_{\lambda}=d_{\mu}$ follows. 
Dipper and James defined Specht modules $S^{\lambda}_K$ for Hecke algebras as irreducible submodules 
of regular modules in \cite{D-J}, and 
improved the theory of representations of Hecke algebras in the series of articles such as 
\cite{D-J,D-J2,D-J3}. 
In particular, they showed in \cite{D-J3} that if $K$ is a field and $\mathcal{H}_{K,r}(q)$ is semisimple, 
then $\hat{S}^{\lambda}_K{\cong_{\mathcal{H}_{K,r}(q)}}S^{\lambda'}_K$ where $\lambda'$ denotes 
the transpose of $\lambda$. 
In this case, $M_{q,\lambda}$ is equivalent to $S^{\lambda}_{\bar{K}}$, thus $\mu=\lambda'$ follows. 
We divide into two cases depending on whether $\lambda$ is self-conjugate or not. \par
\noindent
(Case1) $\lambda{\neq}\lambda'$: \\
In this case, \, $\hat{}$ \, induces an involution on $I_{q,\lambda}{\oplus}I_{q,\lambda'}$. 
Let $\tilde{I}_{q,\lambda}=\big{\{}X+Y{\in}I_{q,\lambda}{\oplus}I_{q,\lambda'}|
(X+Y){\sphat}=X+Y\big{\}}$ for ${\lambda}\neq{\lambda'}$. 
Then we have 
\begin{equation}
\tilde{I}_{q,\lambda}=\big{\{}X+\hat{X}{\in}I_{q,\lambda}{\oplus}I_{q,\lambda'}|
X{\in}I_{q,\lambda}\big{\}}
{\cong}I_{q,\lambda}\big{(}{\cong}\operatorname{Mat}_{d_{\lambda}}(\bar{K})\big{)},
\end{equation}
Thereby, the image of the regular representation of 
$\mathcal{H}^1_{\bar{K},r}(q)$ on $S^{\lambda}_{\bar{K}}$ is isomorphic to 
$I_{q,\lambda}$, so $\operatorname{res}^{\mathcal{H}_{\bar{K},r}(q)}_{\mathcal{H}^1_{\bar{K},r}(q)}
S^{\lambda}_{\bar{K}}$ 
is a simple left $\mathcal{H}^1_{\bar{K},r}(q)$-module. 
If $g{\in}\mathcal{H}_{\bar{K},r}(q)$ satisfies $\hat{g}=g$, then the matrix coefficients with respect to the basis 
$x_1,x_2,\hdots,x_{d_{\lambda}}$ of $S^{\lambda}_{\bar{K}}$ is the same as those with respect to the basis 
$\hat{x}_1,\hat{x}_2,\hdots,\hat{x}_{d_{\lambda}}$ of $\hat{S}^{\lambda}_{\bar{K}}{\cong}
S^{\lambda'}_{\bar{K}}$. 
Hence we have the isomorphism of simple left $\mathcal{H}^1_{\bar{K},r}(q)$-modules as follows. 
\[
\operatorname{res}^{\mathcal{H}_{\bar{K},r}(q)}_{\mathcal{H}^1_{\bar{K},r}(q)}S^{\lambda}_{\bar{K}}{\cong}
\operatorname{res}^{\mathcal{H}_{\bar{K},r}(q)}_{\mathcal{H}^1_{\bar{K},r}(q)}S^{\lambda'}_{\bar{K}}.
\]
(Case2) $\lambda=\lambda'$: \par
\noindent
In this case, \, $\hat{}$ \, induces an involution on 
$I_{q,\lambda}{\cong}\operatorname{Mat}_{d_{\lambda}}(\bar{K})$. 
By Skolem-Noether Theorem, there exists an invertible element $P$ of 
$\operatorname{Mat}_{d_{\lambda}}(\bar{K})$ 
such that $\hat{X}=PXP^{-1}$ for all $X{\in}\operatorname{Mat}_{d_{\lambda}}(\bar{K})$. 
Since eigenvalues of $P$ are ${\pm}1$, We may also assume that 
\[
P=
\begin{bmatrix}
1 &   &   &   &   & \\
  & \ddots &   &   &   & \\
  &   & 1 &   &   & \\
  &   &   & -1 &   & \\
  &   &   &    & \ddots  & \\
  &   &   &    &   & -1
\end{bmatrix}.
\]
Then we may write 
\[
P
\begin{bmatrix}
X_1 & X_2 \\
X_3  & X_4
\end{bmatrix}
P^{-1}=
\begin{bmatrix}
X_1 & -X_2 \\
-X_3  & X_4
\end{bmatrix},
\]
for some submatrices $X_1,X_2,X_3,X_4$. 
Therefore, we obtain 
\[
\big{\{}X{\in}\operatorname{Mat}_{d_{\lambda}}(\bar{K})|\hat{X}=X\big{\}}
=\bigg{\{}X{\in}\operatorname{Mat}_{d_{\lambda}}(\bar{K})|X=
\begin{bmatrix}
X_1 & 0 \\
0  & X_4
\end{bmatrix}\bigg{\}}.
\]
Assuming that $1$ appears $m$ times in $P$. Then 
$\dim_{\bar{K}}\big{\{}X{\in}\operatorname{Mat}_{d_{\lambda}}(\bar{K})|\hat{X}=X\big{\}}
=m^2+(d_{\lambda}-m)^2$. 
We easily see that 
\[
m^2+(d_{\lambda}-m)^2{\geq}\dfrac{d_{\lambda}^2}{2} \quad \text{equality holds iff $m=\dfrac{d_{\lambda}}{2}$}.
\]
Combining this with (3.5) and the fact that 
$\dim_{\bar{K}}\mathcal{H}_{\bar{K},r}^1(q)=\dfrac{1}{2}\dim_{\bar{K}}\mathcal{H}_{\bar{K},r}(q)$ 
(Proposition 3.6), 
we deduce that $m=d_{\lambda}/2$. Thereby, 
\begin{equation}
\big{\{}X{\in}\operatorname{Mat}_{d_{\lambda}}(\bar{K})|\hat{X}=X\big{\}}
{\cong}\operatorname{Mat}_{d_{\lambda}/2}(\bar{K}){\oplus}\operatorname{Mat}_{d_{\lambda}/2}(\bar{K})
\end{equation}
holds. 
This means that 
$\operatorname{res}^{\mathcal{H}_{\bar{K},r}(q)}_{\mathcal{H}_{\bar{K},r}^1(q)}S^{\lambda}_{\bar{K}}$ 
decomposes into two simple left $\mathcal{H}^1_{\bar{K},r}(q)$-modules 
$\bar{S}^{\lambda+}_{\bar{K}}$ and $\bar{S}^{\lambda-}_{\bar{K}}$ which are mutually non-isomorphic. 
Let $\tilde{I}_{q,\lambda}=\big{\{}X{\in}I_{q,\lambda}|\hat{X}=X\big{\}}$ and 
$\tilde{I}_{q,\lambda}^+$ (resp. $\tilde{I}_{q,\lambda}^-$) be the homogeneous component 
corresponding to $\bar{S}^{\lambda+}_{\bar{K}}$ (resp. $\bar{S}^{\lambda-}_{\bar{K}}$) 
for ${\lambda}\neq{\lambda'}$. 
Then (3.6) implies that 
\begin{equation}
\tilde{I}_{q,\lambda}=\tilde{I}_{q,\lambda}^+{\oplus}\tilde{I}_{q,\lambda}^-,
\end{equation}
Summarizing our argument, we conclude that $\mathcal{H}_{\bar{K},r}^1(q)$ is isomorphic to the 
direct sum of minimal two-sided ideals as follows, 
\[
\mathcal{H}_{\bar{K},r}^1(q)=
\big{\{}\bigoplus_{\lambda{\in}\Lambda_r,\lambda>\lambda'}\tilde{I}_{q,\lambda}\big{\}}
\bigoplus\Big{[}\bigoplus_{\lambda{\in}\Lambda_r,\lambda=\lambda'}
\big{\{}\tilde{I}_{q,\lambda}^+{\oplus}\tilde{I}_{q,\lambda}^-\big{\}}\Big{]},
\]
where $<$ denotes the lexicographic order on $\Lambda_r$; $\lambda=(\lambda_1,\lambda_2,\hdots)<
\mu=(\mu_1,\mu_2,\hdots)$ iff $\lambda_k<\mu_k$ for the smallest $k$ such that $\lambda_k{\neq}\mu_k$.  
Consequently we have proved the following result. 
\begin{Theo}
Let $q$ be an indeterminate and $K=\mathbb{Q}(q)$. \\
If $\lambda'{\neq}\lambda$, then $\bar{S}^{\lambda}_{\bar{K}}$ is simple and 
$\bar{S}^{\lambda}_{\bar{K}}{\cong}\bar{S}^{\lambda'}_{\bar{K}}$ 
where $\bar{S}^{\lambda}_{\bar{K}}$ denotes 
$\operatorname{res}^{\mathcal{H}_{\bar{K},r}(q)}_{\mathcal{H}_{\bar{K},r}^1(q)}S^{\lambda}_{\bar{K}}$. \\
If $\lambda'=\lambda$, then $\bar{S}^{\lambda}_{\bar{K}}$ 
decomposes into mutually non-isomorphic two simple left $\mathcal{H}^1_{\bar{K},r}(q)$-modules 
$\bar{S}^{\lambda+}_{\bar{K}}$ and $\bar{S}^{\lambda-}_{\bar{K}}$. \\
The simple left $\mathcal{H}^1_{\bar{K},r}(q)$-modules 
$\bar{S}^{\lambda}_{\bar{K}}, \bar{S}^{\mu+}_{\bar{K}}, \bar{S}^{\mu-}_{\bar{K}} 
(\lambda,\mu{\in}\Lambda_r,\lambda>\lambda',\mu=\mu')$ 
constitute a basic set of simple left $\mathcal{H}^1_{\bar{K},r}(q)$-modules. 
Moreover, $\mathcal{H}_{\bar{K},r}^1(q)$ is a split semisimple 
$\bar{K}$-algebra and its homogeneous decomposition 
is as follows: 
\[
\mathcal{H}_{\bar{K},r}^1(q)=
\big{\{}\bigoplus_{\lambda{\in}\Lambda_r,\lambda>\lambda'}\tilde{I}_{q,\lambda}\big{\}}
\bigoplus\Big{[}\bigoplus_{\lambda{\in}\Lambda_r,\lambda=\lambda'}
\big{\{}\tilde{I}_{q,\lambda}^+{\oplus}\tilde{I}_{q,\lambda}^-\big{\}}\Big{]},
\]
where each $\tilde{I}_{q,\lambda}$ is the homogeneous component corresponding to $\lambda$ which is isomorphic to  
$\operatorname{Mat}_{d_{\lambda}}(\bar{K})$ and each direct sum $\tilde{I}_{q,\lambda}^+{\oplus}\tilde{I}_{q,\lambda}^-$ 
consists of two homogeneous components $\tilde{I}_{q,\lambda}^+$ and $\tilde{I}_{q,\lambda}^-$ both corresponding to 
$\lambda$ which are isomorphic to $\operatorname{Mat}_{d_{\lambda}/2}(\bar{K})$. 
\end{Theo}

\section{Schur-Weyl reciprocity between $U^\sigma_q\big{(}\gl(m,n)\big{)}$ and 
$\mathcal{H}_{K,r}(q)$}
In our previous paper \cite{Mit2}, we defined the $q$-permutation representation and 
established Schur-Weyl reciprocity between the quantum superalgebra $U^\sigma_q\big{(}\gl(m,n)\big{)}$ 
and $\mathcal{H}_{K,r}(q)$. 
Schur-Weyl reciprocity between the general Lie superalgebra $\gl(m,n)$ and the symmetric group 
$\mathfrak{S}_r$ was established in, for example, \cite{B-R,Ser}. 
In this section, we shall review the $q$-permutation representation and the 
vector representation of $U^\sigma_q\big{(}\gl(m,n)\big{)}$. \par
Let $V={\oplus}_{k=1}^{m+n}Kv_k$ be a $\mathbb{Z}_2$-graded $K$-module of rank $m+n$. 
By $\mathbb{Z}_2$-graded, we mean that $V$ is a direct sum of two submodule 
$V_{\overline{0}}=\oplus_{k=1}^{m}Kv_k$ and $V_{\overline{1}}=\oplus_{k=m+1}^{m+n}Kv_k$, and that 
for each homogeneous element the degree map $|{\cdot}|$ 
\begin{equation*}
|v|=
\begin{cases}
0& \text{if $v{\in}V_{\overline{0}}$},\\
1& \text{if $v{\in}V_{\overline{1}}$,}
\end{cases}
\end{equation*}
is given. \par
Let $\pi_r$ be the $q$-permutation representation of $\mathcal{H}_{K,r}(q)$ on the tensor space 
$V^{{\otimes}r}$. 
$\pi_r$ is given by 
$\pi_r(T_i)=\operatorname{Id}^{{\otimes}i-1}{\otimes}T{\otimes}\operatorname{Id}^{{\otimes}r-i-1}$
($i=1,2,\hdots,r-1$) where $T$ is the operator on $V{\otimes}V$ defined by 
\begin{equation}
Tv_k{\otimes}v_l=
\begin{cases}
\dfrac{(-1)^{|v_k|}(q+q^{-1})+q-q^{-1}}{2}v_k{\otimes}v_l& \text{if $k=l$,}\\
(-1)^{|v_k||v_l|}v_l{\otimes}v_k+(q-q^{-1})v_k{\otimes}v_l& \text{if $k<l$,}\\
(-1)^{|v_k||v_l|}v_l{\otimes}v_k& \text{if $k>l$.}
\end{cases}
\end{equation}
and $\operatorname{Id}$ is the identity operator on $V$. 
This representation $\pi_r$ is reduced to the (normal) $q$-permutation representation of 
$\mathcal{H}_{K,r}(q)$ with $n=0$ and to the sign permutation representation 
of $\mathfrak{S}_r$ with $q{\rightarrow}1$. 
Let $T'$ be the operator defined by 
\[
T'=\dfrac{2T-(q-q^{-1})}{q+q^{-1}}
\]
Then $T'$ is determined by 
\begin{equation}
T'v_k{\otimes}v_l=
\begin{cases}
(-1)^{|v_k|}v_k{\otimes}v_l& \text{if $k=l$,}\\
\dfrac{2(-1)^{|v_k||v_l|}}{q+q^{-1}}v_l{\otimes}v_k+\dfrac{q-q^{-1}}{q+q^{-1}}v_k{\otimes}v_l& 
\text{if $k<l$,}\\
\dfrac{2(-1)^{|v_k||v_l|}}{q+q^{-1}}v_l{\otimes}v_k-\dfrac{q-q^{-1}}{q+q^{-1}}v_k{\otimes}v_l& 
\text{if $k>l$.}
\end{cases}
\end{equation}
$\pi_r$ is also given by 
$\pi_r(T'_i)=\operatorname{Id}^{{\otimes}i-1}{\otimes}T'{\otimes}\operatorname{Id}^{{\otimes}r-i-1}$
($i=1,2,\hdots,r-1$)
\par
Next we shall review quantum superalgebras and their vector representations. 
Several definitions of quantum superalgebras appear in, for example, \cite{B-K-K,K-T,Yamane}. 
The much complete definition and detailed observations of quantum superalgebras can be found in 
\cite{Yamane}. In this paper, we obey the manner of definition of $U^\sigma_q\big{(}\gl(m,n)\big{)}$ 
in \cite{B-K-K}. 
The method of construction of superalgebra depends on \cite{Kac} basically. 
Let $\Pi=\{\alpha_i\}_{i{\in}I}$ be a set of simple roots with the index set $I=\{1,\hdots,r\}$. 
We assume that $I$ is a disjoint union of two subsets $I_{\rm even}$ and $I_{\rm odd}$. 
We define a map $p:I{\longrightarrow}\{0,1\}$ to be such that
\begin{equation*}
p(i)=
\begin{cases}
0& \text{if $i{\in}I_{\rm even}$},\\
1& \text{if $i{\in}I_{\rm odd}$.}
\end{cases}
\end{equation*}
Let $P$ be a free $\mathbb{Z}$-module which includes all $\alpha_i{\in}P$($i{\in}I$). 
We assume that a $\mathbb{Q}$-valued symmetric bilinear form on $P$ 
$(\cdot,\cdot) : P{\times}P{\longrightarrow}{\mathbb{Q}}$ is defined and that the 
simple coroots $h_i{\in}P^*$($i{\in}I$) are given as data. The natural pairing 
${\langle}\cdot,\cdot{\rangle} : P^*{\times}P{\longrightarrow}{\mathbb{Z}}$ 
between $P$ and $P^*$ is assumed to satisfy 
\begin{equation*}
{\langle}h_i,\alpha_j{\rangle}=
\begin{cases}
2& \text{if $i=j$ and $i{\in}I_{\rm even}$,}\\
0\quad\text{or}\quad2& \text{if $i=j$ and $i{\in}I_{\rm odd}$,}\\
{\leq}0& \text{if $i{\neq}j$.}
\end{cases}
\end{equation*}
We denote by $\Pi^{\vee}=\{h_i|i{\in}I\}$ the set of all coroots. 
Furthermore, for each $i{\in}I$ we assume that there exists a nonzero integer $\ell_i$ such that 
$\ell_i{\langle}h_i,\lambda{\rangle}=(\alpha_i,\lambda)$ for every $\lambda{\in}P$.
Then we immediately have the Cartan matrix $A=[{\langle}h_i,\alpha_j{\rangle}]_{ij}$ 
is symmetrizable because $\ell_i{\langle}h_i,\alpha_j{\rangle}=(\alpha_i,\alpha_j)
=(\alpha_j,\alpha_i)=\ell_j{\langle}h_j,\alpha_i{\rangle}$. We mention that 
the symmetrized matrix is $A^{\rm sym}=\operatorname{diag}(\ell_1,\hdots,\ell_r)A=
[(\alpha_i,\alpha_j)]_{ij}$. Let $\mathfrak{h}=P^*{\otimes}_{\mathbb{Z}}\mathbb{Q}$. 
Then $\Phi=(\mathfrak{h},\Pi^{\vee},\Pi)$ 
is said to be a fundamental root data associated to $A$. 
Let $\g=\g(\Phi)$ be the contragredient Lie superalgebra obtained from $\Phi$ and $p$. 
The quantized enveloping algebra $U_q^\sigma(\g)$ is the unital associative algebra over 
$K=\mathbb{Q}(q)$ with generators $q^h (h{\in}P^*),e_i,f_i (i{\in}I)$ and an additional element 
$\sigma$ which satisfy the following defining relations:
\renewcommand{\theenumi}{\arabic{enumi}}
\renewcommand{\labelenumi}{(Q\theenumi)}
\begin{enumerate}
\item
$q^h=1$ \quad for $h=0$, 
\item
$q^{h_1}q^{h_2}=q^{h_1+h_2}$ \quad for $h_1,h_2{\in}P^*$, 
\item
$q^he_i=q^{{\langle}h,\alpha_j{\rangle}}e_iq^h$ \quad for $h{\in}P^*$ and $i{\in}I$,
\item
$q^hf_i=q^{-{\langle}h,\alpha_j{\rangle}}f_iq^h$ \quad for $h{\in}P^*$ and $i{\in}I$,
\item
$[e_i,f_j]=\delta_{ij}\dfrac{q^{\ell_ih_i}-q^{-\ell_ih_i}}{q^{\ell_i}-q^{-\ell_i}}$ 
\quad for $i,j{\in}I$,
\item
$\sigma^2=1$,
\item
$q^h{\sigma}={\sigma}q^h$ \quad for $h{\in}P^*$,
\item
$e_i{\sigma}=(-1)^{p(i)}{\sigma}e_i$ \quad for $i{\in}I$,
\item
$f_i{\sigma}=(-1)^{p(i)}{\sigma}f_i$ \quad for $i{\in}I$, 
\end{enumerate}
\renewcommand{\theenumi}{\arabic{enumi}}
\renewcommand{\labelenumi}{(\theenumi)}
where $[e_i,f_j]$ means the supercommutator 
\begin{equation*}
[e_i,f_j]=e_if_j-(-1)^{p(i)p(j)}f_je_i.
\end{equation*}
We assume further conditions:
\renewcommand{\theenumi}{\arabic{enumi}}
\renewcommand{\labelenumi}{(Q\theenumi)}
\begin{enumerate}
\item[(Q10)]
If $a{\in}\sum_{i{\in}I}U_q(\mathfrak{n}_+)e_iU_q(\mathfrak{n}_+)$ satisfies 
$f_ia{\in}U_q(\mathfrak{n}_+)f_i$ for all $i{\in}I$, then $a=0$, 
\item[(Q11)]
If $a{\in}\sum_{i{\in}I}U_q(\mathfrak{n}_-)f_iU_q(\mathfrak{n}_-)$ satisfies 
$e_ia{\in}U_q(\mathfrak{n}_-)e_i$ for all $i{\in}I$, then $a=0$, 
\end{enumerate}
\renewcommand{\theenumi}{\arabic{enumi}}
\renewcommand{\labelenumi}{(\theenumi)}
where $U_q(\mathfrak{n}_+)$ (respectively $U_q(\mathfrak{n}_-)$) is the subalgebra of 
$U_q^\sigma(\g)$ generated by $\{e_i|i{\in}I\}$ (respectively $\{f_i|i{\in}I\}$). 
$U_q^\sigma(\g)$ is a Hopf algebra whose comultiplication $\triangle_{\sigma}$, 
counit $\varepsilon_{\sigma}$, antipode $S_{\sigma}$ are as follows. 
\begin{equation*}
\begin{split}
&\triangle_{\sigma}(\sigma)=\sigma{\otimes}\sigma, \\
&\triangle_{\sigma}(q^h)=q^h{\otimes}q^h \quad \text{for $h{\in}P^*$}, \\
&\triangle_{\sigma}(e_i)=e_i{\otimes}q^{-\ell_ih_i}+\sigma^{p(i)}{\otimes}e_i \quad \text{for $i{\in}I$}, \\
&\triangle_{\sigma}(f_i)=f_i{\otimes}1+\sigma^{p(i)}q^{\ell_ih_i}{\otimes}f_i \quad \text{for $i{\in}I$}, \\
&\varepsilon_{\sigma}(\sigma)=\varepsilon_{\sigma}(q^h)=1 
\quad \text{for $h{\in}P^*$}, \quad 
\varepsilon_{\sigma}(e_i)=\varepsilon_{\sigma}(f_i)=0 \quad \text{for $i{\in}I$}, \\
&S_{\sigma}(\sigma)=\sigma, \quad S_{\sigma}(q^{{\pm}h})=q^{{\mp}h} 
\quad \text{for $h{\in}P^*$}, \\
&S_{\sigma}(e_i)=-\sigma^{p(i)}e_iq^{\ell_ih_i},{\quad}S_{\sigma}(f_i)=-\sigma^{p(i)}q^{-\ell_ih_i}f_i 
\quad \text{for $i{\in}I$}. 
\end{split}
\end{equation*}
The quantized enveloping algebra $U^\sigma_q\big{(}\gl(m,n)\big{)}$ 
is obtained from the fundamental root data as follows. 
\begin{itemize}
\item
$I=I_{\rm even}{\cup}I_{\rm odd}$ is defined by 
$I_{\rm even}=\{1,2,\hdots,m-1,m+1,\hdots,m+n-1\}$ and $I_{\rm odd}=\{m\}$,
\item
$P=\oplus_{b{\in}B}\mathbb{Z}\epsilon_b$, where 
$B=B_+{\cup}B_-$ with $B_+=\{1,\hdots,m\}$ and $B_-=\{m+1,\hdots,m+n\}$,
\item
$(\cdot,\cdot):P{\times}P{\longrightarrow}\mathbb{Q}$ is 
the symmetric bilinear form on $P$ defined by
\[
(\epsilon_a,\epsilon_{a'})=
\begin{cases}
1& \text{if $a=a'{\in}B_+$,}\\
-1& \text{if $a=a'{\in}B_-$,}\\
0& \text{otherwise,}
\end{cases}
\]
\item
$\Pi=\{\alpha_i|i{\in}I\}$ is defined by $\alpha_i=\epsilon_i-\epsilon_{i+1}$,
\item
$\Pi^{\vee}=\{h_i|i{\in}I\}$ is uniquely determined by the formula 
$\ell_i{\langle}h_i,\lambda{\rangle}=(\alpha_i,\lambda)$ for any $\lambda{\in}P$, \\
where 
\[
\ell_i=
\begin{cases}
1& \text{if $i=1,\hdots,m$,}\\
-1& \text{if $i=m+1,\hdots,m+n-1$.}
\end{cases}
\]
\end{itemize}\par
The vector representation ($\rho,V$) of $U^\sigma_q\big{(}\gl(m,n)\big{)}$ on the 
$\mathbb{Z}_2$-graded vector space $V=V_{\overline{0}}{\oplus}V_{\overline{1}}$ (recall that 
$V_{\overline{0}}=\oplus_{i=1}^{m}Kv_i,V_{\overline{1}}=\oplus_{i=m+1}^{m+n}Kv_i$) is defined by 
(see \cite{B-K-K}) 
\begin{equation}
\begin{split}
&\rho(\sigma)v_j=(-1)^{|v_j|}v_j \quad \text{for $j=1,\hdots,m+n$}, \\
&\rho(q^h)v_j=q^{\epsilon_j(h)}v_j \quad \text{for $h{\in}P^*,j=1,\hdots,m+n$}, \\
&\rho(e_j)v_{j+1}=v_j \quad \text{for $j=1,\hdots,m+n-1$}, \\
&\rho(f_j)v_j=v_{j+1} \quad \text{for $j=1,\hdots,m+n-1$}, \\
&\text{otherwise $0$.}
\end{split}
\end{equation}
The vector representation ${\rho}_r$ on $V^{{\otimes}r}$ is given 
by ${\rho}_r(x)=\rho^{{\otimes}r}\circ\triangle^{(r-1)}(x)$ where 
$\triangle^{(1)}=\triangle_\sigma$ and 
$\triangle^{(k)}=(\triangle_\sigma{\otimes}\operatorname{Id}^{{\otimes}k-1})\triangle^{(k-1)}$ 
inductively. One can readily see that $\rho_r$ is of the following form: 
\begin{equation}
\begin{split}
&{\rho}_r(\sigma)=\rho(\sigma)^{{\otimes}r}, \\
&{\rho}_r(q^h)=\rho(q^h)^{{\otimes}r} \quad \text{for $h{\in}P^*$}, \\
&{\rho}_r(e_i)=\sum_{k=1}^N\rho(\sigma^{p(i)})^{{\otimes}k-1}{\otimes}\rho(e_i) 
{\otimes}\rho(q^{-\ell_ih_i})^{{\otimes}r-k} \quad 
\text{for $i{\in}I$}, \\
&{\rho}_r(f_i)=\sum_{k=1}^r\rho(\sigma^{p(i)}q^{\ell_ih_i})^{{\otimes}k-1}{\otimes}\rho(f_i)
{\otimes}\operatorname{Id}^{{\otimes}r-k} \quad 
\text{for $i{\in}I$}.
\end{split}
\end{equation}
Our precedent works in \cite{Mit2} are Theorem 4.1 and Theorem 4.2. 
\begin{Theo}[Schur-Weyl reciprocity \cite{Mit2}]
Let $\mathcal{A}_q=\pi_r\big{(}\mathcal{H}_{K,r}(q)\big{)}$ and 
$\mathcal{B}_q=\rho_r\big{(}U^\sigma_q\big{(}\gl(m,n)\big{)}\big{)}$. 
Then we have $\operatorname{End}_{\mathcal{B}_q}V^{{\otimes}r}
=\mathcal{A}_q$ and $\operatorname{End}_{\mathcal{A}_q}V^{{\otimes}r}
=\mathcal{B}_q$. 
\end{Theo}
Let $\bar{U}^\sigma_q\big{(}\gl(m,n)\big{)}=U^\sigma_q\big{(}\gl(m,n)\big{)}
{\otimes_{K}}\bar{K}$, 
$\bar{\mathcal{A}}_q=\mathcal{A}_q{\otimes_{K}}\bar{K}$ and 
$\bar{\mathcal{B}}_q=\mathcal{B}_q{\otimes_{K}}\bar{K}$. 
Then, 
$\pi_r\big{(}\mathcal{H}_{\bar{K},r}(q)\big{)}=\bar{\mathcal{A}}_q$ and 
$\rho_r\big{(}\bar{U}^\sigma_q\big{(}\gl(m,n)\big{)}\big{)}=\bar{\mathcal{B}}_q$ 
as $\bar{K}$-algebras of operators on 
$\bar{V}^{{\otimes}r}=(V{\otimes_{K}}\bar{K})^{{\otimes}r}$. 
We notice that $\operatorname{End}_{\bar{\mathcal{B}}_q}\bar{V}^{{\otimes}r}
=\bar{\mathcal{A}}_q$ and $\operatorname{End}_{\bar{\mathcal{A}}_q}\bar{V}^{{\otimes}r}
=\bar{\mathcal{B}}_q$ hold. \par
Let $H(m,n;r)=\{\lambda=(\lambda_1,\lambda_2,\hdots){\in}\Lambda_r|\lambda_j{\leq}n 
\text{ if } j>m\}$. Diagrams of elements of $H(m,n;r)$ are exactly those contained 
in the $(m,n)$-hooks. Then the following holds. 
\begin{Theo}[Decomposition of the tensor space \cite{Mit2}]
$\bar{\mathcal{A}}_q=\bigoplus_{\lambda{\in}H(m,n;r)}\bar{\mathcal{A}}_{q,\lambda}$ 
where each $\bar{\mathcal{A}}_{q,\lambda}=\pi_r(I_{q,\lambda})({\cong}I_{q,\lambda})$ is the image of 
the homogeneous component $I_{q,\lambda}$ corresponding to $\lambda{\in}H(m,n;r)$ as in (3.4). Moreover, we have 
$\bar{V}^{{\otimes}r}=\bigoplus_{\lambda{\in}H(m,n;r)}H_{\lambda}{\otimes}V_{\lambda}$ 
where $H_{\lambda}$'s are mutually non-isomorphic simple left 
$\mathcal{H}_{\bar{K},r}(q)$-modules indexed by the elements of $H(m,n;r)$, and 
$V_{\lambda}$'s are mutually non-isomorphic simple left 
$\bar{U}^\sigma_q\big{(}\gl(m,n)\big{)}$-modules indexed by the elements of $H(m,n;r)$. 
\end{Theo}

\section{Schur-Weyl reciprocity for $\mathcal{H}_{K,r}^1(q)$ in case of $m=n$}
Schur-Weyl reciprocity for the alternating group $\mathfrak{A}_r$ has been researched by Regev 
in \cite{Regev}. 
The paper \cite{Regev} showed that if $\dim V_{\overline{0}}=\dim V_{\overline{1}}$ under the base 
field $\mathbb{C}$, then the centralizer algebra $\operatorname{End}_{\mathfrak{A}_r}V^{{\otimes}r}$ 
has remarkable property; $\mathbb{Z}_2$-crossed product for 
$\operatorname{End}_{\mathfrak{S}_r}V^{{\otimes}r}$. 
In this paper, we establish a $q$-analogue extension of Regev's result. 
We also show that Schur-Weyl reciprocity is valid even if the base field is $\mathbb{Q}(q)$. \par
We set $K=\mathbb{Q}(q)$ in succession, and denote 
$\mathcal{C}_q=\pi_r\big{(}\mathcal{H}_{K,r}^1(q)\big{)}$. 
Let us consider the relation between 
$\operatorname{End}_{\mathcal{A}_q}V^{{\otimes}r}$ 
and $\operatorname{End}_{\mathcal{C}_q}V^{{\otimes}r}$. 
Since $\mathcal{H}_{K,r}^1(q){\subsetneq}\mathcal{H}_{K,r}(q)$, 
we immediately have 
\[
\mathcal{C}_q{\subseteq}\mathcal{A}_q \quad 
\text{and} \quad 
\operatorname{End}_{\mathcal{A}_q}V^{{\otimes}r}{\subseteq}
\operatorname{End}_{\mathcal{C}_q}V^{{\otimes}r}.
\]
Recall that $V$ is an $m+n$-dimensional $\mathbb{Z}_2$-graded vector space over $K$. 
In this section, we analyze the structure of $\operatorname{End}_{\mathcal{C}_q}V^{{\otimes}r}$ 
in case of $m=n$. We will consider the general case in the next section. \par
Assume that $m=n$. Recall that 
\begin{equation*}
\mathcal{B}_q=\{f{\in}\operatorname{End}_{K}V^{{\otimes}r}
|\pi_r(T'_i)f=f\pi_r(T'_i) \text{ for $i=1,2,\hdots,r-1$}\} 
=\operatorname{End}_{\mathcal{A}_q}V^{{\otimes}r}.
\end{equation*}
We define $\mathcal{B}_q^{\dagger},\mathcal{C}_q$ to be the subspaces of 
$\operatorname{End}_{K}V^{{\otimes}r}$ as follows. 
\begin{equation*}
\begin{split}
\mathcal{B}_q^{\dagger}&=\{f{\in}\operatorname{End}_{K}V^{{\otimes}r}
|\pi_r(T'_i)f=-f\pi_r(T'_i) \text{ for $i=1,2,\hdots,r-1$}\}, \\
\mathcal{D}_q&=\{f{\in}\operatorname{End}_{K}V^{{\otimes}r}
|\pi_r(T'_1T'_{i+1})f=f\pi_r(T'_1T'_{i+1}) \text{ for $i=1,2,\hdots,r-2$}\} 
=\operatorname{End}_{\mathcal{C}_q}V^{{\otimes}r}.
\end{split}
\end{equation*}
\begin{Lem}
$\mathcal{D}_q=\mathcal{B}_q{\bigoplus}\mathcal{B}_q^{\dagger}$.
\end{Lem}
\begin{proof}
It is clear that the sum is direct. 
$\supseteq$ is also obvious. 
We notice that if $f{\in}\mathcal{D}_q$, then $\pi_r(T'_j)f\pi_r(T'_j)=\pi_r(T'_1)f\pi_r(T'_1)$ for 
$j=2,3,\hdots,r-1$ from the definition of $\mathcal{D}_q$. \par
In general, we may write
\[
f=\dfrac{1}{2}\big{(}f+\pi_r(T'_1)f\pi_r(T'_1)\big{)}+\dfrac{1}{2}\big{(}f-\pi_r(T'_1)f\pi_r(T'_1)\big{)}.
\]
If $f{\in}\mathcal{D}_q$, then one can readily see 
\begin{equation*}
\begin{split}
\pi_r(T'_i)\dfrac{1}{2}\big{(}f+\pi_r(T'_1)f\pi_r(T'_1)\big{)}
&=\dfrac{1}{2}\big{(}\pi_r(T'_i)f\pi_r(T'_i)^2+\pi_r(T'_i)\pi_r(T'_1)f\pi_r(T'_1)\big{)}\\
&=\dfrac{1}{2}\big{(}\pi_r(T'_1)f\pi_r(T'_1)\pi_r(T'_i)+f\pi_r(T'_i)\pi_r(T'_1)^2\big{)}\\
&=\dfrac{1}{2}\big{(}\pi_r(T'_1)f\pi_r(T'_1)+f\big{)}\pi_r(T'_i).
\end{split}
\end{equation*}
Hence we have $2^{-1}\big{(}f+\pi_r(T'_1)f\pi_r(T'_1)\big{)}{\in}\mathcal{B}_q$. 
In the same fashion, we also get $2^{-1}\big{(}f-\pi_r(T'_1)f\pi_r(T'_1)\big{)}{\in}\mathcal{B}_q^{\dagger}$. 
Thus we have proved the reverse inclusion $\subseteq$. 
\end{proof}
We notice that if $V$ is $\mathbb{Z}_2$-graded, then $\operatorname{End}_{K}V$ is also 
$\mathbb{Z}_2$-graded. Namely, 
\[
\operatorname{End}_{K}V=(\operatorname{End}_{K}V)_{\overline{0}}
{\oplus}(\operatorname{End}_{K}V)_{\overline{1}},
\]
where 
\begin{equation*}
\begin{split}
(\operatorname{End}_{K}V)_{\overline{0}}
&=\big{\{}f{\in}\operatorname{End}_{K}V|f(V_{\overline{i}}){\subseteq}V_{\overline{i}},
i{\in}\{0,1\}\big{\}}, \\
(\operatorname{End}_{K}V)_{\overline{1}}
&=\big{\{}f{\in}\operatorname{End}_{K}V|f(V_{\overline{i}}){\subseteq}V_{\overline{i+1}},
i{\in}\{0,1\}\big{\}}.
\end{split}
\end{equation*}
Let $\varphi{\in}\operatorname{End}_{K}V$ be given by 
$\varphi(v_i)=v_{2m-i+1}$. 
Obviously, $\varphi{\in}(\operatorname{End}_{K}V)_{\overline{1}}$. 
Let $\varphi^{{\otimes}r}$ be the tensor product of $\varphi$. 
In general, for homogeneous elements $f_1,f_2,\hdots,f_r$ of $\operatorname{End}_{K}V$, 
$f_1{\otimes}f_2{\otimes}\cdots{\otimes}f_r{\in}\operatorname{End}_{K}V^{{\otimes}r}$ 
is given by 
\[
f_1{\otimes}\cdots{\otimes}f_r(u_1{\otimes}\cdots{\otimes}u_r)
=(-1)^{\sum_{i=2}^r(|f_i|\sum_{j=1}^{i-1}|u_j|)}
f_1(u_1){\otimes}\cdots{\otimes}f_r(u_r),
\]
where $u_1,\hdots,u_r$ are homogeneous elements of $V$. 
Hence we have 
\begin{equation}
\varphi^{{\otimes}r}(u_1{\otimes}\cdots{\otimes}u_r)
=(-1)^{\sum_{i=2}^r(\sum_{j=1}^{i-1}|u_j|)}
\varphi(u_1){\otimes}\cdots{\otimes}\varphi(u_r).
\end{equation}
\begin{Lem}
$\varphi^{{\otimes}r}$ is an isomorphism of $K$-vector space $V^{{\otimes}r}$ 
which satisfies following properties:  
\begin{equation}
(\varphi^{{\otimes}r})^2=(-1)^{r(r-1)/2}\operatorname{I}
\end{equation}
and
\begin{equation}
\pi_r(T'_i)\varphi^{{\otimes}r}=-\varphi^{{\otimes}r}\pi_r(T'_i)\quad\text{for $i=1,2,\hdots,r-1$},
\end{equation}
where $\operatorname{I}$ is the identity operator on $V^{{\otimes}r}$.
\end{Lem}
\begin{proof}
For homogeneous elements $u_1,\hdots,u_r$ of $V$, one can readily see that 
\begin{equation*}
\begin{split}
{(\varphi^{{\otimes}r})^2}(u_1{\otimes}\cdots{\otimes}u_r)
&=(-1)^{\sum_{i=2}^r(\sum_{j=1}^{i-1}|u_j|)}
\varphi(u_1){\otimes}\cdots{\otimes}\varphi(u_r) \\
&=(-1)^{\sum_{i=2}^r(\sum_{j=1}^{i-1}(|u_j|+|\varphi(u_j)|)}
u_1{\otimes}\cdots{\otimes}u_r \\
&=(-1)^{\sum_{i=2}^r\sum_{j=1}^{i-1}1}
u_1{\otimes}\cdots{\otimes}u_r \\
&=(-1)^{(r-1)r/2}
u_1{\otimes}\cdots{\otimes}u_r.
\end{split}
\end{equation*}
Hence $(\varphi^{{\otimes}r})^2=(-1)^{r(r-1)/2}\operatorname{I}$, and $\varphi^{{\otimes}r}$ 
has the inverse $(\varphi^{{\otimes}r})^{-1}=(-1)^{(r-1)r/2}\varphi^{{\otimes}r}$. \par
To prove the last statement, 
we check three cases of the definition of $\pi_r$ which has appeared in (4.2). 
We notice that it suffice to prove only the case $r=2$. \\
Case1 : $k=l$ 
\begin{equation*}
\begin{split}
{\varphi^{{\otimes}2}}T'(v_k{\otimes}v_l)&=(-1)^{|v_k|}\varphi^{{\otimes}2}(v_k{\otimes}v_k) \\
&=v_{2m-k+1}{\otimes}v_{2m-k+1}
\end{split}
\end{equation*}
\begin{equation*}
\begin{split}
T'{\varphi^{{\otimes}2}}(v_k{\otimes}v_l)&=(-1)^{|v_k|}T'(v_{2m-k+1}{\otimes}v_{2m-k+1}) \\
&=(-1)^{|v_k|}(-1)^{|v_{2m-k+1}|}v_{2m-k+1}{\otimes}v_{2m-k+1} \\
&=-v_{2m-k+1}{\otimes}v_{2m-k+1}
\end{split}
\end{equation*}
Case2 : $k<l$
\begin{equation*}
\begin{split}
{\varphi^{{\otimes}2}}T'(v_k{\otimes}v_l)&=\dfrac{2(-1)^{|v_k||v_l|}}{q+q^{-1}}{\varphi^{{\otimes}2}}(v_l{\otimes}v_k)
+\dfrac{q-q^{-1}}{q+q^{-1}}{\varphi^{{\otimes}2}}(v_k{\otimes}v_l) \\
&=\dfrac{2(-1)^{|v_k||v_l|+|v_l|}}{q+q^{-1}}(v_{2m-l+1}{\otimes}v_{2m-k+1})
+\dfrac{q-q^{-1}}{q+q^{-1}}(-1)^{|v_k|}(v_{2m-k+1}{\otimes}v_{2m-l+1})
\end{split}
\end{equation*}
\begin{equation*}
\begin{split}
T'{\varphi^{{\otimes}2}}(v_k{\otimes}v_l)&=(-1)^{|v_k|}T'(v_{2m-k+1}{\otimes}v_{2m-l+1}) \\
&=\dfrac{2(-1)^{|v_k|+|v_{2m-k+1}||v_{2m-l+1}|}}{q+q^{-1}}(v_{2m-l+1}{\otimes}v_{2m-k+1})
-\dfrac{q-q^{-1}}{q+q^{-1}}(-1)^{|v_k|}(v_{2m-k+1}{\otimes}v_{2m-l+1})
\end{split}
\end{equation*}
If $|v_k|=|v_l|=0$, then $(-1)^{|v_k||v_l|+|v_l|}=1$ and 
$(-1)^{|v_k|+|v_{2m-k+1}||v_{2m-l+1}|}=-1$.\\
If $|v_k|=0,|v_l|=1$, then $(-1)^{|v_k||v_l|+|v_l|}=-1$ and 
$(-1)^{|v_k|+|v_{2m-k+1}||v_{2m-l+1}|}=1$.\\
If $|v_k|=1,|v_l|=0$, then $(-1)^{|v_k||v_l|+|v_l|}=1$ and 
$(-1)^{|v_k|+|v_{2m-k+1}||v_{2m-l+1}|}=-1$.\\
If $|v_k|=|v_l|=1$, then $(-1)^{|v_k||v_l|+|v_l|}=1$ and 
$(-1)^{|v_k|+|v_{2m-k+1}||v_{2m-l+1}|}=-1$.\\
After all, ${\varphi^{{\otimes}2}}T'(v_k{\otimes}v_l)=-T'{\varphi^{{\otimes}2}}(v_k{\otimes}v_l)$ holds. \\
Case3 : $k>l$
\begin{equation*}
\begin{split}
{\varphi^{{\otimes}2}}T'(v_k{\otimes}v_l)&=\dfrac{2(-1)^{|v_k||v_l|}}{q+q^{-1}}{\varphi^{{\otimes}2}}(v_l{\otimes}v_k)
-\dfrac{q-q^{-1}}{q+q^{-1}}{\varphi^{{\otimes}2}}(v_k{\otimes}v_l) \\
&=\dfrac{2(-1)^{|v_k||v_l|+|v_l|}}{q+q^{-1}}(v_{2m-l+1}{\otimes}v_{2m-k+1})
-\dfrac{q-q^{-1}}{q+q^{-1}}(-1)^{|v_k|}(v_{2m-k+1}{\otimes}v_{2m-l+1})
\end{split}
\end{equation*}
\begin{equation*}
\begin{split}
T'{\varphi^{{\otimes}2}}(v_k{\otimes}v_l)&=(-1)^{|v_k|}T'(v_{2m-k+1}{\otimes}v_{2m-l+1}) \\
&=\dfrac{2(-1)^{|v_k|+|v_{2m-k+1}||v_{2m-l+1}|}}{q+q^{-1}}(v_{2m-l+1}{\otimes}v_{2m-k+1})
+\dfrac{q-q^{-1}}{q+q^{-1}}(-1)^{|v_k|}(v_{2m-k+1}{\otimes}v_{2m-l+1})
\end{split}
\end{equation*}
In the same manner as case2, we have ${\varphi^{{\otimes}2}}T'(v_k{\otimes}v_l)=-T'{\varphi^{{\otimes}2}}(v_k{\otimes}v_l)$. 
\end{proof}
Let ${\Phi}$ be the endomorphism of the vector space $\operatorname{End}_{K}V^{{\otimes}r}$ 
which is given by: 
\[
{\Phi}(f)={\varphi^{{\otimes}r}}f \quad 
\text{for $f{\in}\operatorname{End}_{K}V^{{\otimes}r}$}.
\]
Then we have the following. 
\begin{Prop}
$\Phi$ is an automorphism of the vector space $\operatorname{End}_{K}V^{{\otimes}r}$. 
The restriction of $\Phi$ to $\mathcal{D}_q$ gives an 
automorphism of $\mathcal{D}_q$ which satisfies the following properties:
\begin{equation*}
{\Phi}(\mathcal{B}_q)=\mathcal{B}_q^{\dagger}, \quad {\Phi}(\mathcal{B}_q^{\dagger})=\mathcal{B}_q.
\end{equation*}
Especially, we have 
$\dim_{K}\mathcal{B}_q=\dim_{K}\mathcal{B}_q^{\dagger}$ and 
$\mathcal{D}_q=\mathcal{B}_q{\oplus}\Phi(\mathcal{B}_q)=\mathcal{B}_q{\oplus}\varphi^{{\otimes}r}\mathcal{B}_q$. 
\end{Prop}
\begin{proof}
From (5.2), it immediately follows that ${\Phi^2}f=(-1)^{r(r-1)/2}f$ for 
$f{\in}\operatorname{End}_{K}V^{{\otimes}r}$. 
Hence $\Phi$ is an isomorphism of $\operatorname{End}_{K}V^{{\otimes}r}$. 
Let $f{\in}\mathcal{B}_q$. From (5.3), it follows that  
\begin{equation*}
\begin{split}
{\pi_r}(T'_i){\Phi}(f)&={\pi_r}(T'_i){\varphi^{{\otimes}r}}f \\
&=-{\varphi^{{\otimes}r}}{\pi_r}(T'_i)f \\
&=-{\varphi^{{\otimes}r}}f{\pi_r}(T'_i) \\
&=-{\Phi}(f){\pi_r}(T'_i).
\end{split}
\end{equation*}
Hence ${\Phi}(f){\in}\mathcal{B}_q^{\dagger}$. 
In the same way, we also obtain ${\Phi}(f){\in}\mathcal{B}_q$ for $f{\in}\mathcal{B}_q^{\dagger}$. 
Combining with Lemma 5.1, we can conclude that $\Phi$ defines an 
automorphism of $\mathcal{D}_q=\mathcal{B}_q{\oplus}\mathcal{B}_q^{\dagger}$ 
satisfying ${\Phi}(\mathcal{B}_q)=\mathcal{B}_q^{\dagger}$ and 
${\Phi}(\mathcal{B}_q^{\dagger})=\mathcal{B}_q$. 
\end{proof}
Let $\omega$ be the endomorphism of the algebra $\operatorname{End}_{K}V^{{\otimes}r}$ 
which is given by: 
\begin{equation*}
{\omega}(f)=(-1)^{r(r-1)/2}{\varphi^{{\otimes}r}}f{\varphi^{{\otimes}r}} \quad 
\text{for $f{\in}\operatorname{End}_{K}V^{{\otimes}r}$}.
\end{equation*}
Indeed, one can easily see that the following. 
\begin{equation*}
\begin{split}
{\omega}(fg)&=(-1)^{r(r-1)/2}{\varphi^{{\otimes}r}}fg{\varphi^{{\otimes}r}} \\
&={\varphi^{{\otimes}r}}f{\varphi^{{\otimes}r}}{\varphi^{{\otimes}r}}g{\varphi^{{\otimes}r}} \\
&=\big{(}(-1)^{r(r-1)/2}{\varphi^{{\otimes}r}}f{\varphi^{{\otimes}r}}\big{)}
\big{(}(-1)^{r(r-1)/2}{\varphi^{{\otimes}r}}g{\varphi^{{\otimes}r}}\big{)} \\
&={\omega}(f){\omega}(g).
\end{split}
\end{equation*}
So, $\omega$ is algebraic. 
\begin{Prop}
$\omega$ is an automorphism of order $2$ of the algebra $\mathcal{B}_q$.
\end{Prop}
\begin{proof}
If $f{\in}\mathcal{B}_q$, then it follows from (5.3) that 
\begin{equation*}
\begin{split}
{\pi_r(T'_i)}{\omega}(f)&=(-1)^{r(r-1)/2}{\pi_r(T'_i)}{\varphi^{{\otimes}r}}f{\varphi^{{\otimes}r}} \\
&=(-1)^{r(r-1)/2}{\pi_r(T'_i)}{\varphi^{{\otimes}r}}f{\varphi^{{\otimes}r}} \\
&=(-1)^{r(r-1)/2}{\varphi^{{\otimes}r}}f{\varphi^{{\otimes}r}}{\pi_r(T'_i)} \\
&={\omega}(f){\pi_r(T'_i)}
\end{split}
\end{equation*}
Therefore ${\omega}(f){\in}\mathcal{B}_q$. 
From (5.2), it follows that 
\begin{equation*}
\begin{split}
{\omega^2}(f)&=(-1)^{r(r-1)/2}{\omega}({\varphi^{{\otimes}r}}f{\varphi^{{\otimes}r}}) \\
&=(-1)^{r(r-1)}{(\varphi^{{\otimes}r})^2}f{(\varphi^{{\otimes}r})^2} \\
&=f
\end{split}
\end{equation*}
Therefore $\omega$ is an automorphism of $\mathcal{B}_q$ of order $2$. 
\end{proof}
Let $H={\langle}1,\omega{\rangle}$. Then $H$ is a subgroup of $\operatorname{Aut}(\mathcal{B}_q)$. 
$H$ is naturally isomorphic to $\mathbb{Z}_2={\langle}1,-1{\rangle}$ as (multiplicative) groups. 
We define two maps 
\begin{equation}
\psi_1:\mathbb{Z}_2{\longrightarrow}\operatorname{Aut}(\mathcal{B}_q),\quad 
\psi_1(1)=1,\psi_1(-1)=\omega,
\end{equation}
and 
\begin{equation}
\alpha_1:\mathbb{Z}_2{\times}\mathbb{Z}_2{\longrightarrow}(\mathcal{B}_q)^{\times},\quad 
\alpha_1(\sigma,\tau)=
\begin{cases}
1& \text{if $\sigma=1$ or $\tau=1$,}\\
(-1)^{r(r-1)/2}& \text{otherwise.}
\end{cases}
\end{equation}
Then we have the following lemma. 
\begin{Lem}
$\psi_1$ and $\alpha_1$ satisfy (2.1)--(2.3).
\end{Lem}
\begin{proof}
(2.3) is trivial. Since $\alpha(\sigma,\tau)={\pm}1$, (2.1) may be reduced to 
${^\sigma}({^\tau}a)={^{({\sigma}{\tau})}}a$, so holds obviously. 
If $\sigma_1=1$, then both sides of (2.2) equal $\alpha(\sigma_2,\sigma_3)$. 
In the same manner, if $\sigma_2=1$ (respectively $\sigma_3=1$), 
then those equal $\alpha(\sigma_1,\sigma_3)$ (respectively $\alpha(\sigma_1,\sigma_2)$). 
If $\sigma_1=\sigma_2=\sigma_3=-1$, then both sides of (2.2) become $-1$. 
Thus (2.2) holds. 
\end{proof}
Hence, we obtain a crossed product ${\mathcal{B}_q}_{\alpha_1}^{\psi_1}[\mathbb{Z}_2]$ from 
the crossed system $(\mathcal{B}_q,\mathbb{Z}_2,\psi_1,\alpha_1)$, 
where $\psi_1$ and $\alpha_1$ are given in (5.4) and (5.5) respectively. 
\begin{Theo}
If $m=n$, then $\mathcal{D}_q$ is isomorphic to the $\mathbb{Z}_2$-crossed product 
${\mathcal{B}_q}_{\alpha_1}^{\psi_1}[\mathbb{Z}_2]$ as $K$-algebras. 
\end{Theo}
\begin{proof}
Because $\mathcal{B}_q^{\dagger}=\varphi^{{\otimes}r}\mathcal{B}_q=\mathcal{B}_q\varphi^{{\otimes}r}$ 
holds from (5.3) and Proposition 5.3, 
\[
\mathcal{D}_q=\mathcal{B}_q1{\oplus}\mathcal{B}_q\varphi^{{\otimes}r}
\]
follows. 
Since both ${\mathcal{B}_q}_{\alpha_1}^{\psi_1}[\mathbb{Z}_2]$ and $\mathcal{D}_q$ are free left 
$\mathcal{B}_q$-modules, we may define an isomorphism of $\mathcal{B}_q$-modules: 
\begin{equation*}
{\iota_1}:{\mathcal{B}_q}_{\alpha_1}^{\psi_1}[\mathbb{Z}_2]{\longrightarrow}\mathcal{D}_q,\quad 
{\iota_1}(u_1)=1,{\iota_1}(u_{-1})=\varphi^{{\otimes}r}.
\end{equation*}
From (2.4) and (5.4) and (5.5), one can deduce the multiplication law for 
${\mathcal{B}_q}_{\alpha_1}^{\psi_1}[\mathbb{Z}_2]$ as follows. 
\begin{equation*}
(a_1u_{\sigma})(a_2u_{\tau})=
\begin{cases}
a_1a_2u_1& \text{if $\sigma=1$ and $\tau=1$,}\\
a_1a_2u_{-1}& \text{if $\sigma=1$ and $\tau=-1$,}\\
(-1)^{r(r-1)/2}a_1\varphi^{{\otimes}r}a_2\varphi^{{\otimes}r}u_{-1}& \text{if $\sigma=-1$ and $\tau=1$,}\\
a_1\varphi^{{\otimes}r}a_2\varphi^{{\otimes}r}u_1& \text{if $\sigma=-1$ and $\tau=-1$,}\\
\end{cases}
\end{equation*}
where $a_1,a_2{\in}\mathcal{B}_q$. Therefore, we get four formulas in $\mathcal{D}_q$, 
\begin{eqnarray*}
(a_11)(a_21)&=&a_1a_21,\\
(a_11)(a_2\varphi^{{\otimes}r})&=&a_1a_2\varphi^{{\otimes}r},\\
(a_1\varphi^{{\otimes}r})(a_21)&=&(-1)^{r(r-1)/2}a_1\varphi^{{\otimes}r}a_2
\varphi^{{\otimes}r}\varphi^{{\otimes}r},\\
(a_1\varphi^{{\otimes}r})(a_2\varphi^{{\otimes}r})&=&a_1\varphi^{{\otimes}r}a_2
\varphi^{{\otimes}r}1,
\end{eqnarray*}
which derive the conclusion that $f$ is an isomorphism of $K$-algebras. 
\end{proof}
For a commutative domain $R$ and a subalgebra $A$ of $\operatorname{Mat}(m,R)$, 
we set $\Tilde{A}=\{X{\in}\operatorname{Mat}(m,R)|XY=YX \text{ for all }Y{\in}A\}$. 
Let $R_1=\mathbb{Q}[q,q^{-1}]$ be the algebra of Laurent polynomials. 
Since matrix elements of $\rho_r$ are in $R_1$ by (4.3) and (4.4), 
we may define $\mathcal{D}'_q$ to be the subalgebra of 
$\operatorname{Mat}\big{(}(2m)^r,{R_1}\big{)}$ generated by the set 
$\{\rho_r(\sigma),\rho_r(q^h),\rho_r(e_i),\rho_r(f_i),\varphi^{{\otimes}r}|h{\in}P^*,i{\in}I\}$. 
Similarly by (4.2), we may also define $\mathcal{C}'_q$ to be the one generated by 
$\{I_{(2m)^r},(q+q^{-1})^2\pi_r(T'_1T'_i)|i=2,\hdots,r-1\}$ where $I_{(2m)^r}$ is the identity matrix. 
Now we shall complete Schur-Weyl reciprocity for the $q$-analogue of the alternating group, namely, 
$\mathcal{C}_q=\operatorname{End}_{\mathcal{D}_q}V^{{\otimes}r}$. \par
$\mathcal{D}_q=\Tilde{\mathcal{C}}_q$ is by definition. 
We shall show $\mathcal{C}_q=\Tilde{\mathcal{D}}_q$. 
The specialization to a nonzero complex number $t$ is a ring homomorphism 
$\varphi_t : R_1{\longrightarrow}\mathbb{C}$ with the condition $\varphi_t(q)=t$. 
$\mathbb{C}$ becomes $(\mathbb{C},R_1)$-bimodule, with $R_1$ acting from the right via 
$\varphi_t$. 
If $t$ is a transcendental number, we can extend the specialization from $R_1$ to its quotient field 
$K$, namely 
$\varphi_t : K{\longrightarrow}\mathbb{C}$. 
Applying the specialization $\varphi_t$, we obtain the specialized 
algebras $\mathcal{C}_t=\mathbb{C}{\otimes_{R_1}}\mathcal{C}'_q$ and 
$\mathcal{D}_t=\mathbb{C}{\otimes_{R_1}}\mathcal{D}_q$ which are subalgebras of 
$\operatorname{Mat}\big{(}(2m)^r,\mathbb{C}\big{)}$. 
$\mathcal{C}_t$ and $\mathcal{D}_t$ act on the specialized vector space 
$\mathbb{C}{\otimes}_KV^{{\otimes}r}$ in obvious ways. 
\begin{Prop}
$\mathcal{C}_q=\Tilde{\mathcal{D}}_q$. 
\end{Prop}
\begin{proof}
Obviously $\mathcal{C}'_q{\subseteq}\Tilde{\mathcal{D}}'_q$, so 
$\rank_{R_1}\mathcal{C}'_q{\leq}\rank_{R_1}\Tilde{\mathcal{D}}'_q$ holds. \par
Since $R_1$ is a principal ideal domain, the submodules $\mathcal{C}'_q$ and $\mathcal{D}'_q$ 
of the free $R_1$-module $\operatorname{Mat}\big{(}(2m)^r,R_1\big{)}$ are also free. 
Assume $N=\rank_{R_1}\mathcal{C}'_q$. Let $X_i(q)$($i=1,\hdots,N$) be a basis of 
$\mathcal{C}'_q$ and $x_i^{k,l}(q){\in}R_1$ the $(k,l)$-entry of $X_i(q)$. 
Then we immediately have that the specialized elements 
$X_i(t)=\big{(}x_i^{k,l}(t)\big{)}$($i=1,\hdots,N$) generate $\mathcal{C}_t$ and 
$\dim_\mathbb{C}\mathcal{C}_t{\leq}\rank_{R_1}\mathcal{C}'_q$. 
Because $X_i(q)$ are linearly independent, 
$\sum_{i=1}^N\alpha_i(q)X_i(q)=0$ for $\alpha_1(q),\hdots,\alpha_N(q){\in}R_1$ 
implies $\alpha_1(q)=\hdots=\alpha_N(q)=0$. 
But the equation $\sum_{i=1}^N\alpha_iX_i(t)=0$ for $\alpha_1,\hdots,\alpha_N{\in}\mathbb{C}$ 
may admit a nonzero $\alpha_i$. 
Consider the system of linear equations:  
\begin{equation*}
\begin{bmatrix}
x_1^{1,1}(t) & x_2^{1,1}(t) & \cdots & x_N^{1,1}(t) \\
x_1^{1,2}(t) & x_2^{1,2}(t) & \cdots & x_N^{1,2}(t) \\
&\cdots\cdots\cdots \\
x_1^{(2m)^r,(2m)^r-1}(t) & x_2^{(2m)^r,(2m)^r-1}(t) & \cdots & x_N^{(2m)^r,(2m)^r-1}(t) \\
x_1^{(2m)^r,(2m)^r}(t) & x_2^{(2m)^r,(2m)^r}(t) & \cdots & x_N^{(2m)^r,(2m)^r}(t)
\end{bmatrix}
\times
\begin{bmatrix}
\alpha_1 \\
\alpha_2 \\
\cdots \\
\alpha_N
\end{bmatrix}
=
\begin{bmatrix}
0 \\
0 \\
\cdots \\
0 \\
0
\end{bmatrix}.
\end{equation*}
Reducing to a common denominator for each row, we may assume that each $x_i^{k,l}(t){\in}\mathbb{Q}[t]$. 
Suppose there exists a non-trivial solution. 
Then the rank of above matrix is less than $N$, 
which implies that every $N$th minor determinant equals $0$. 
This is impossible if $t$ is transcendental. 
Hence $\dim_\mathbb{C}\mathcal{C}_t=\rank_{R_1}\mathcal{C}'_q$ for a transcendental $t$. 
In the same manner, $\dim_\mathbb{C}\mathcal{D}_t=\rank_{R_1}\mathcal{D}'_q$ 
if $t$ is a transcendental number. 
This argument is valid even if $R_1$ is replaced by $K$, thus we also obtain 
$\dim_\mathbb{C}\mathcal{C}_t=\dim_K\mathcal{C}_q$, 
$\dim_\mathbb{C}\mathcal{D}_t=\dim_K\mathcal{D}_q$ 
for a transcendental number $t$. \par
Next we shall show that $\rank_{R_1}\Tilde{\mathcal{D}}'_q=
\dim_\mathbb{C}\Tilde{\mathcal{D}}_t$. 
Assume that $Y(q)=\big{(}y^{k,l}(q)\big{)}{\in}\Tilde{\mathcal{D}}'_q$. 
Then $Y(q)$ commutes with all 
$\rho_r(\sigma),\rho_r(q^h),\rho_r(e_i),\rho_r(f_i),\varphi^{{\otimes}r}(h{\in}P^*,i{\in}I)$. 
Hence the matrix elements $y^{k,l}(q)(k,l=1,\hdots,(2m)^r$ determine polynomial 
equations with coefficients in $R_1$. 
We notice that $\rank_{R_1}\Tilde{\mathcal{D}}'_q{\leq}\dim_\mathbb{C}\Tilde{\mathcal{D}}_t$ in 
general; 
besides the elements of specialized algebra $\mathbb{C}{\otimes_{R_1}}\Tilde{\mathcal{D}}'_q$, 
possibly $\Tilde{\mathcal{D}}_t$ contains another element which commutes with specializations of 
$\rho_r(\sigma),\rho_r(q^h),\rho_r(e_i),\rho_r(f_i),\varphi^{{\otimes}r}(h{\in}P^*,i{\in}I)$. 
Assume the case. Reducing to a common denominator, 
we may deduce that $t$ is a common solution of certain polynomial equations with coefficients 
in $\mathbb{Q}$. This is impossible if $t$ is transcendental. 
Thus $\rank_{R_1}\Tilde{\mathcal{C}}'_q=\dim_\mathbb{C}\Tilde{\mathcal{C}}_t$ for a transcendental $t$. 
In the same manner, we also have 
$\rank_{R_1}\Tilde{\mathcal{D}}'_q=\dim_\mathbb{C}\Tilde{\mathcal{D}}_t$ for a transcendental $t$. 
Replacing $R_1$ by $K$, we obtain 
$\dim_K\Tilde{\mathcal{C}}_q=\dim_\mathbb{C}\Tilde{\mathcal{C}}_t$ and 
$\dim_K\Tilde{\mathcal{D}}_q=\dim_\mathbb{C}\Tilde{\mathcal{D}}_t$. 
Considering the specialization $\varphi_1$, we obtain the following: 
\begin{equation*}
\dim_{\mathbb{C}}\mathcal{C}_1{\leq}
\dim_{\mathbb{C}}\mathcal{C}_t=
\rank_{R_1}\mathcal{C}'_q{\leq}
\rank_{R_1}\Tilde{\mathcal{D}}'_q=
\dim_{\mathbb{C}}\Tilde{\mathcal{D}}_t{\leq}
\dim_{\mathbb{C}}\Tilde{\mathcal{D}}_1.
\end{equation*}
In \cite{Regev}, Regev has shown $\mathcal{D}_1=\Tilde{\mathcal{C}}_1$. 
Applying double centralizer theorem, we have $\mathcal{C}_1=\Tilde{\mathcal{D}}_1$. 
Therefore $\rank_{R_1}\mathcal{C}'_q=\rank_{R_1}\Tilde{\mathcal{D}}'_q$ and 
$\dim_K\mathcal{C}_q=\dim_K\Tilde{\mathcal{D}}_q$, which imply $\mathcal{C}_q=\Tilde{\mathcal{D}}_q$. 
\end{proof}
From Proposition 5.7 we immediately obtain Schur-Weyl reciprocity for $\mathcal{H}_{K,r}^1(q)$ 
as follows. 
\begin{Theo}
$\operatorname{End}_{\mathcal{C}_q}V^{{\otimes}r}
=\mathcal{D}_q$ and $\operatorname{End}_{\mathcal{D}_q}V^{{\otimes}r}
=\mathcal{C}_q$ hold. 
\end{Theo}

\section{Schur-Weyl reciprocity for $\mathcal{H}_{\bar{K},r}^1(q)$ in the general case}
Henceforth, we do not assume $m{\neq}n$ and consider the general case. 
In this case, results can be obtained in the same way as $q=1$ case which is used in 
\cite{Regev}. 
Recall Theorem 4.2:
\begin{equation}
\bar{\mathcal{A}}_q=\bigoplus_{\lambda{\in}H(m,n;r)}\bar{\mathcal{A}}_{q,\lambda},
\end{equation}
where each $\bar{\mathcal{A}}_{q,\lambda}=\pi_r(I_{q,\lambda})({\cong}I_{q,\lambda})$ is the image of 
the homogeneous component $I_{q,\lambda}$ of $\mathcal{H}_{\bar{K},r}(q)$ 
corresponding to $\lambda{\in}H(m,n;r)$. 
In the same manner as \cite{Regev}, we define 
$H_0(m,n;r)=\{\lambda{\in}\Lambda_r|\lambda,\lambda'{\in}H(m,n;r)\}$ and 
$H_1(m,n;r)=H(m,n;r){\backslash}H_0(m,n;r)$. 
Then we obtain the following from (6.1).  
\begin{equation}
\bar{\mathcal{A}}_q=\Big{[}\bigoplus_{\lambda{\in}H_0(m,n;r),\lambda>\lambda'}
\big{\{}\bar{\mathcal{A}}_{q,\lambda}{\oplus}\bar{\mathcal{A}}_{q,\lambda'}\big{\}}\Big{]}
\bigoplus\big{\{}\bigoplus_{\lambda{\in}H_0(m,n;r),\lambda=\lambda'}\bar{\mathcal{A}}_{q,\lambda}\big{\}}
\bigoplus\big{\{}\bigoplus_{\lambda{\in}H_1(m,n;r)}\bar{\mathcal{A}}_{q,\lambda}\big{\}}.
\end{equation}
Let $\bar{\mathcal{C}}_q=\mathcal{C}_q{\otimes_{\mathbb{Q}(q)}}\bar{K}
=\pi_r\big{(}\mathcal{H}_{\bar{K},r}^1(q)\big{)}$. 
From Theorem 3.10 we have:
\[
\mathcal{H}_{\bar{K},r}^1(q)=
\big{\{}\bigoplus_{\lambda{\in}\Lambda_r,\lambda>\lambda'}\tilde{I}_{q,\lambda}\big{\}}
\bigoplus\Big{[}\bigoplus_{\lambda{\in}\Lambda_r,\lambda=\lambda'}
\big{\{}\tilde{I}_{q,\lambda}^+{\oplus}\tilde{I}_{q,\lambda}^-\big{\}}\Big{]}.
\]
So we immediately obtain:
\begin{equation}
\bar{\mathcal{C}}_q=
\big{\{}\bigoplus_{\lambda{\in}H_0(m,n;r),\lambda>\lambda'}\pi_r(\tilde{I}_{q,\lambda})\big{\}}
\bigoplus\Big{[}\bigoplus_{\lambda{\in}H_0(m,n;r),\lambda=\lambda'}
\big{\{}\pi_r(\tilde{I}_{q,\lambda}^+){\oplus}\pi_r(\tilde{I}_{q,\lambda}^-)\big{\}}\Big{]}
\bigoplus\big{\{}\bigoplus_{\lambda{\in}H_1(m,n;r)}\pi_r(\tilde{I}_{q,\lambda})\big{\}}.
\end{equation}
If $\lambda{\in}H_0(m,n;r)$ and $\lambda>\lambda'$, then $\pi_r(I_{q,\lambda}){\cong}I_{q,\lambda}$ 
and $\pi_r(I_{q,\lambda'}){\cong}I_{q,\lambda'}$. Therefore 
$\pi_r(I_{q,\lambda}{\oplus}I_{q,\lambda'}){\cong}I_{q,\lambda}{\oplus}I_{q,\lambda'}$, which implies 
$\pi_r(\tilde{I}_{q,\lambda}){\cong}\pi_r(I_{q,\lambda}){\cong}I_{q,\lambda}$ because 
$I_{q,\lambda}{\cong}\tilde{I}_{q,\lambda}{\subseteq}I_{q,\lambda}{\oplus}I_{q,\lambda'}$. \par
If $\lambda{\in}H_0(m,n;r)$ and $\lambda=\lambda'$, then 
$\pi_r(\tilde{I}_{q,\lambda}^+){\cong}\tilde{I}_{q,\lambda}^+$ and 
$\pi_r(\tilde{I}_{q,\lambda}^-){\cong}\tilde{I}_{q,\lambda}^-$ hold because 
$\pi_r(\tilde{I}_{q,\lambda}){\cong}\tilde{I}_{q,\lambda}$ and 
$\tilde{I}_{q,\lambda}^{\pm}{\subseteq}\tilde{I}_{q,\lambda}$. \par
If $\lambda{\in}H_1(m,n;r)$ then $\pi_r(I_{q,\lambda}){\cong}I_{q,\lambda}$ and 
$\pi_r(I_{q,\lambda'})=\{0\}$. Since 
$\tilde{I}_{q,\lambda}=\big{\{}X+\hat{X}{\in}I_{q,\lambda}{\oplus}I_{q,\lambda'}|
X{\in}I_{q,\lambda}\big{\}}$, we have $\pi_r(X+\hat{X})=\pi_r(X)$ for 
$X+\hat{X}{\in}\tilde{I}_{q,\lambda}$. Thereby $\pi_r(\tilde{I}_{q,\lambda})=\pi_r(I_{q,\lambda})
{\cong}I_{q,\lambda}$ follows. \par
Let us define two-sided ideals of $\bar{\mathcal{A}}_q$ and those of $\bar{\mathcal{C}}_q$ 
as follows. 
\begin{eqnarray*}
\bar{\mathcal{A}}_q{\supseteq}\bar{\mathcal{A}}_q^0&=&
\Big{[}\bigoplus_{\lambda{\in}H_0(m,n;r),\lambda>\lambda'}
\big{\{}\bar{\mathcal{A}}_{q,\lambda}{\oplus}\bar{\mathcal{A}}_{q,\lambda'}\big{\}}\Big{]}
{\bigoplus}\big{\{}\bigoplus_{\lambda{\in}H_0(m,n;r),\lambda=\lambda'}
\bar{\mathcal{A}}_{q,\lambda}\big{\}}, \\
\bar{\mathcal{A}}_q{\supseteq}\bar{\mathcal{A}}_q^1&=&\bigoplus_{\lambda{\in}H_1(m,n;r)}
\bar{\mathcal{A}}_{q,\lambda}, \\
\bar{\mathcal{C}}_q{\supseteq}\bar{\mathcal{C}}_q^0&=&
\big{\{}\bigoplus_{\lambda{\in}H_0(m,n;r),\lambda>\lambda'}
\pi_r(\tilde{I}_{q,\lambda})\big{\}}{\bigoplus}\Big{[}\bigoplus_{\lambda{\in}H_0(m,n;r),\lambda=\lambda'}
\big{\{}\pi_r(\tilde{I}_{q,\lambda}^+){\oplus}\pi_r(\tilde{I}_{q,\lambda}^-)\big{\}}\Big{]} \\
\bar{\mathcal{C}}_q{\supseteq}\bar{\mathcal{C}}_q^1&=&\bigoplus_{\lambda{\in}H_1(m,n;r)}
\pi_r(\tilde{I}_{q,\lambda})
\end{eqnarray*}
Comparing dimensions of components between $\bar{\mathcal{A}}_q$ and 
$\bar{\mathcal{C}}_q$, we obtain the following theorem. 
\begin{Theo}
$\bar{\mathcal{A}}_q$ and $\bar{\mathcal{C}}_q$ have direct sum decompositions 
$\bar{\mathcal{A}}_q=\bar{\mathcal{A}}_q^0{\oplus}\bar{\mathcal{A}}_q^1$ and 
$\bar{\mathcal{C}}_q=\bar{\mathcal{C}}_q^0{\oplus}\bar{\mathcal{C}}_q^1$ respectively, 
which are satisfy the following relations. 
\begin{enumerate}
\item
$\bar{\mathcal{C}}_q^0{\subseteq}\bar{\mathcal{A}}_q^0$ and 
$\dim_{\bar{K}}\bar{\mathcal{A}}_q^0
=2\dim_{\bar{K}}\bar{\mathcal{C}}_q^0$
\item
$\bar{\mathcal{C}}_q^1=\bar{\mathcal{A}}_q^1$. Especially 
$\dim_{\bar{K}}\bar{\mathcal{A}}_q^1
=\dim_{\bar{K}}\bar{\mathcal{C}}_q^1$. 
\end{enumerate}
\end{Theo}
As a special case, when $m=n$, we readily see that $H_1(m,n;r)=\emptyset$. Hence 
both $\bar{\mathcal{A}}_q^1$ and $\bar{\mathcal{C}}_q^1$ are zero and 
$\dim_{\bar{K}}\bar{\mathcal{A}}_q=2\dim_{\bar{K}}\bar{\mathcal{C}}_q$. 
Moreover, if $n=0$, namely non-super case, 
then we readily see that if $m^2<r$, then all $\lambda{\in}\Lambda_r$ are contained in 
$H_1(m,0;r)$ and $\bar{\mathcal{A}}_q^0=\bar{\mathcal{C}}_q^0=0$. 
The Corollary below is immediately obtained from this argument and (6.2), (6.3). 
\begin{Cor}
Let $n=0$. If $m^2<r$, then $\bar{\mathcal{A}}_q=\bar{\mathcal{C}}_q$ and 
$\operatorname{End}_{\bar{\mathcal{C}}_q}\bar{V}^{{\otimes}r}=
\operatorname{End}_{\bar{\mathcal{A}}_q}\bar{V}^{{\otimes}r}$. 
\end{Cor}
We denote by $W_0$ and $W_1$ direct summands of $\bar{V}^{{\otimes}r}$ which are defined as follows. 
\begin{eqnarray*}
W_0&=&\Big{[}\bigoplus_{\lambda{\in}H_0(m,n;r),\lambda>\lambda'}
\big{\{}(H_{\lambda}{\otimes}V_{\lambda}){\oplus}(H_{\lambda'}{\otimes}V_{\lambda'})\big{\}}\Big{]}
{\bigoplus}\big{\{}\bigoplus_{\lambda{\in}H_0(m,n;r),\lambda=\lambda'}
(H_{\lambda}{\otimes}V_{\lambda})\big{\}} \\
W_1&=&\bigoplus_{\lambda{\in}H_1(m,n;r)}(H_{\lambda}{\otimes}V_{\lambda})
\end{eqnarray*}
We notice that $\operatorname{End}_{\bar{\mathcal{A}}_q}\bar{V}^{{\otimes}r}=
\operatorname{End}_{\bar{\mathcal{A}}_q}W_0{\oplus}\operatorname{End}_{\bar{\mathcal{A}}_q}W_1$ and that 
$\operatorname{End}_{\bar{\mathcal{C}}_q}\bar{V}^{{\otimes}r}=
\operatorname{End}_{\bar{\mathcal{C}}_q}W_0{\oplus}\operatorname{End}_{\bar{\mathcal{C}}_q}W_1$. 
From Theorem 6.1, we immediately have the following. 
\begin{Cor}
$\operatorname{End}_{\bar{\mathcal{C}}_q}W_0{\supseteq}\operatorname{End}_{\bar{\mathcal{A}}_q}W_0$ and 
$\operatorname{End}_{\bar{\mathcal{C}}_q}W_1=\operatorname{End}_{\bar{\mathcal{A}}_q}W_1$.
\end{Cor}
\begin{proof}
Since $\bar{\mathcal{A}}_q^1W_0=0$ and $\bar{\mathcal{A}}_q^0W_1=0$, 
$\operatorname{End}_{\bar{\mathcal{A}}_q}W_0=
\operatorname{End}_{\bar{\mathcal{A}}_q^0}W_0$ and 
$\operatorname{End}_{\bar{\mathcal{A}}_q}W_1=
\operatorname{End}_{\bar{\mathcal{A}}_q^1}W_1$ hold. 
Similarly, we have $\operatorname{End}_{\bar{\mathcal{C}}_q}W_0=
\operatorname{End}_{\bar{\mathcal{C}}_q^1}W_0$ and $\operatorname{End}_{\bar{\mathcal{C}}_q}W_1=
\operatorname{End}_{\bar{\mathcal{C}}_q^1}W_1$. 
Hence we immediately obtain from Theorem 6.1(1) that  
$\operatorname{End}_{\bar{\mathcal{C}}_q}W_0=\operatorname{End}_{\bar{\mathcal{C}}_q^0}W_0
{\supseteq}\operatorname{End}_{\bar{\mathcal{A}}_q^0}W_0=\operatorname{End}_{\bar{\mathcal{A}}_q}W_0$ 
and from Theorem 6.1(2) that 
$\operatorname{End}_{\bar{\mathcal{C}}_q}W_1=\operatorname{End}_{\bar{\mathcal{C}}_q^1}W_1=
\operatorname{End}_{\bar{\mathcal{A}}_q^1}W_1=\operatorname{End}_{\bar{\mathcal{A}}_q}W_1$ also. 
\end{proof}


\begin{thebibliography}{99}
\newcounter{rmnum}
\setcounter{rmnum}{3}
\bibitem{C-R1}
	C. W. Curtis and I. Reiner, $\lq\lq$Methods of Representation Theory", Vol.1, 
		John Wily \& Sons, 1981.
\bibitem{B-K-K}
	G. Benkart, S. Kang and M. Kashiwara, Crystal bases for the quantum superalgebra
	$U_q(\gl(m,n))$, 
	\it J. Amer. Math. Soc. \bf 13(2) \rm (2000), 295--331.
\bibitem{B-R}
	A. Berele and A. Regev, Hook Young diagrams with applications to combinatorics and to 
	representations of Lie superalgebras, 
	\it Adv. Math. \bf 64 \rm (1987), 118--175.
\bibitem{D-J}
	R. Dipper and G.D. James, Representations of Hecke algebras of general linear groups, 
	\it Proc. London Math. Soc. \bf 52 \rm (1986), 20--52.
\bibitem{D-J2}
	R. Dipper and G.D. James, Blocks and idempotents of Hecke algebras of general linear groups, 
	\it Proc. London Math. Soc. \bf 54 \rm (1987), 57--82.
\bibitem{D-J3}
	R. Dipper and G.D. James, The $q$-Schur algebra, 
	\it Proc. London Math. Soc. \bf 59 \rm (1989), 23--50.
\bibitem{Kac}
	V. G. Kac, Lie superalgebras, 
		\it Adv. in Math. \bf 26 \rm (1977), 8--96.
\bibitem{K-T}
	S. M. Khoroshkin and V. N. Tolstoy, Universal $R$-matrix for quantized (super)algebras, 
	\it Commun. Math. Phys. \bf 141 \rm (1991), 599-617.
\bibitem{Mit}
	H. Mitsuhashi, The $q$-analogue of the alternating group and its representations, 
	\it J. Algebra \bf 240 \rm (2001), 535--558.
\bibitem{Mit2}
	H. Mitsuhashi, Schur-Weyl reciprocity between the quantum superalgebra 
	and the Iwahori-Hecke algebra, 
	\it preprint \rm.
\bibitem{N-O}
	C. Nastasescu and F. V. Oystaeyen, Methods of Graded Rings, Lecture Notes in Math., 
		No. 1836, Springer-Verlag, Berlin, 2004.
\bibitem{Regev}
	A. Regev, Double centralizing theorems for the alternating groups, 
	\it J. Algebra \bf 250 \rm (2002), 335--352.
\bibitem{Ser}
	A. N. Sergeev, The tensor algebra of the identity representation as a module over 
	the Lie superalgebras $\mathfrak{G}\mathfrak{l}(n,m)$ and $Q(n)$, 
	\it Math. USSR Sbornik. \bf 51 \rm (1985), No.2, 419--427.
\bibitem{Yamane}
	H. Yamane, Quantized enveloping algebras associated with simple Lie superalgebras and 
	their universal $R$-matrices, 
	\it Publ. RIMS. \bf 30 \rm (1994), 15--87.
\end{thebibliography}
\end{document}